%% file: ex_article_SIOPT_v4.tex
\documentclass[final,hidelinks,onefignum,onetabnum]{siamart220329}

\allowdisplaybreaks

\input{ex_shared}
\expandafter\patchcmd\csname\string\algorithmic\endcsname%
      {\labelwidth 0.5em}{\labelwidth0pt\labelsep0pt}{}{}

\makeatletter
\newcommand{\leqnomode}{\tagsleft@true}
\newcommand{\reqnomode}{\tagsleft@false}
\makeatother
\usepackage{cancel}

\ifpdf
\hypersetup{
  pdftitle={A Modified Proximal Bundle Method Under A Frank-Wolfe Perspective},
}
\fi

\begin{document}

\usetikzlibrary{arrows.meta, intersections, quotes, angles, quotes}
\setlength\fboxsep{0pt}
\renewcommand\algorithmicindent{0.5em}

\maketitle

\begin{abstract}
The proximal bundle method (PBM) is a fundamental and computationally effective algorithm for solving optimization problems with nonsmooth components. In this paper, we conduct a theoretical investigation of a modified proximal bundle method, which we call the Modified Proximal Bundle with Fixed Absolute Accuracy (MPB-FA). MPB-FA modifies the \textit{serious-step test} of the PBM \textit{serious-step test}. In MPB-FA, it is based on an absolute accuracy criterion, and the accuracy is fixed over iterations, while the standard PBM uses a relative accuracy in the serious-step test, which changes with iterations. Also, similarly to multiple PBM analyses, the proximal parameter in MPB-FA is also fixed over iterations, while it is permitted to change in the standard PBM. These modifications allow us to {build the first link between}  a proximal bundle method {and} a Frank-Wolfe algorithm on the Moreau envelope of the dual problem. In light of this correspondence, we first extend the linear convergence of Kelley's method on {the sum of a smooth strongly convex function and a} convex piecewise linear function from the positive homogeneous to the general case. Building on this result, we propose a novel complexity analysis of MPB-FA {when the objective is the sum of a smooth and a piecewise function} and derive an {$\cO(\epsilon^{-8/9})$} iteration complexity, improving upon the best known $\cO(\epsilon^{-2})$ guarantee on a related variant of PBM. It is worth-noting that the best known complexity bound for the classical PBM in the general case is $\cO(\epsilon^{-3})$ {and $\cO(\epsilon^{-2})$ when the proximal parameter is fixed}. Our approach also reveals new insights on bundle management {and empirical behavior of the proximal bundle methods}.
\end{abstract}

\begin{keywords}
iteration-complexity, proximal bundle method, optimal complexity bound
\end{keywords}

\begin{MSCcodes}
90C25, 90C30, 90C60, 90C46  
\end{MSCcodes}

\section{Introduction}

\subsection{Problem statement}
Central to numerous machine learning algorithms \cite{Bundle_teo10a, JMLR:v7:shalev-shwartz06a,Bundle_Methods_ML_2007} and optimization algorithms \cite{bertsimas2022robust, SunConejo2021, BenTal2009, Kuchlbauer2022} lies the convex composite problem 
\begin{align}
\underset{x \in \R^n}{\min} \,\, {h(x):=} g(x) + f(x), \label{problem: primal}
\end{align}
{where $g:\R^n\rightarrow\R$ is a convex smooth function, i.e., $g$ has Lipschitz continuous gradients and $f:\R^n\rightarrow\R$ is a convex piecewise linear function represented as
\begin{align}
f(x) = \max_{(v,b) \in \cV}\{v^\top x + b\} \label{eq: f}
\end{align}
for some finite set $\cV\in\R^{n+1}$. Throughout the paper, we assume \eqref{problem: primal} has an optimal solution $x^*$ and we want to find an $\epsilon$-optimal solution $x$ such that $h(x)-h(x^*)\le\epsilon$.} 


Kelley's method \cite{Kelley1960} is a fundamental algorithm for solving such non-smooth problems (see \Cref{alg: Kelley's algorithm}). 
However, the complexity bound of Kelley's method grows exponentially with the dimension of the problem \cite{Nesterov2004}. In \cite{zhou2018limited}, it is demonstrated that imposing more stringent hypotheses on \eqref{problem: primal}, specifically requiring \( f \) to be \textit{positive homogeneous} and \( g \) to be strongly convex, results in a linear convergence rate for Kelley's method. Although it is feasible to relax the additional hypothesis on \( f \) (as shown in \Cref{section: Linear convergence of Kelley's algorithm}), the strong convexity of \( g \) remains a key requirement in the proof.

In order to accommodate more general $g$, we turn to the Proximal Bundle Method (PBM) \cite{Bonnans2006}. PBM introduces two additional key features to Kelley's method: a quadratic proximal term in the objective and a \textit{serious-step test} to update the proximal center only if the algorithm has made sufficient progress. 
These two key features stabilize the convergence of PBM. In this paper, we attempt to provide a new perspective for PBM via the lens of a Frank-Wolfe algorithm. In particular,  
we modify PBM in two aspects to make a precise dual correspondence to a Frank-Wolfe algorithm: 1) we modify {the \textit{serious-step test}} from using a relative accuracy criterion in PBM to using a fixed absolute accuracy criterion; 2) we keep the proximal parameter constant during iterations, while PBM allows it to change. We call the resulting algorithm the Modified Proximal Bundle with Fixed Absolute Accuracy (MPB-FA). Although MPB-FA is more restrictive than PBM, it inherits two defining features of PBM, namely, it collects a bundle of dual information for the non-smooth objective function and uses a quadratic proximal term for primal stabilization \cite{Kiwiel1990, Bonnans2006}. 

In this paper, we theoretically investigate MPB-FA for the composite problem \eqref{problem: primal}. We give an interpretation of the MPB-FA's \textit{null steps} as dual to the Fully Corrective Frank Wolfe (FCFW) algorithm applied to the Moreau envelope of the conjugate of $g$. This new perspective yields a convergence rate of {${\cO({\epsilon^{-8/9}}(\log{\epsilon}^{-1})^{2/9})}$} for MPB-FA, a significant improvement over the previous rate of {$\cO(\epsilon^{-2})$}. 

\subsection{Related Work}


\textbf{Kelley's method and bundle methods}.
Kelley’s method \cite{Kelley1960} is a foundational algorithm for minimizing non-smooth convex functions. It iteratively refines a piecewise linear lower approximation of the true function. The PBM \cite{Kiwiel1990,deOliveira2014}, the trust region bundle method \cite{Schramm1992}, and the {level bundle} method \cite{Lemarechal1995, Lan2015} {can be viewed as stabilization} of Kelley's method. Recently, Diaz and Grimmer \cite{Diaz2023optimal} showed that their parallel PBM is on par with a state-of-the-art sub-gradient solver \cite{ShalevShwartz2011} for binary SVMs and \cite{atenas:hal-03738298} adapted PBM to convex multistage programming problems.

\textbf{The proximal bundle method}. In 2000, Kiwiel \cite{Kiwiel2000} established the first convergence rate for the PBM, proving that an $\epsilon$-minimizer can be found in $\cO(\epsilon^{-3})$ iterations. This result relies on mild assumptions that accommodate varying proximal parameters and convex constraints. Updating the proximal parameter dynamically has been shown to be crucial for obtaining good numerical performance (see \cite{Bonnans2006} and references therein). The study of proximal bundle type algorithms with fixed proximal parameter $\rho$, especially during sequences of \textit{null steps}, in \cite{Diaz2023optimal, liang2021proximal, LiangMonteiro2023} simplified the analysis and enabled new optimal optimal convergence rates compared to previous analysis. 
\cite{Du2017}, then \cite{Diaz2023optimal} improved the rate of \cite{Kiwiel2000} and provided convergence rates for all combinations of smoothness and strong convexity. That is, $\cO(\epsilon^{-2})$ when the function is only Lipschitz continuous, $\cO(\epsilon^{-1})$ when the function is smooth or strongly convex. When the function is smooth and strongly convex, the convergence rate becomes $\cO(\log\epsilon^{-1})$. In \cite{liang2021proximal, LiangMonteiro2023} the authors studied the composite problem \eqref{problem: primal} and proposed a new type of {\textit{serious step test}} that we have adapted in this paper. Their analysis provides an $\cO(\epsilon^{-2}\log\epsilon^{-1})$ convergence rate when $f$ is only Lipschitz continuous, $\cO(\epsilon^{-1}\log\epsilon^{-1})$ when $f$ is supposed to be smooth. When $g$ is strongly convex the corresponding rates become,  respectively, $\cO(\epsilon^{-1}(\log \epsilon^{-1})^2)$ and $\cO((\log\epsilon^{-1})^2)$. Note that $g$ {may be extended-real-valued}, in contrast to both our setup and the setups in \cite{Du2017, Diaz2023optimal}, and thus could include constraints. 

\textbf{Frank-Wolfe}.
The Frank–Wolfe algorithm \cite{FrankWolfe1956} is a first-order optimization method for constrained convex problems. It does not require projection onto the feasible region; rather, it progresses toward an extreme point determined by a linear minimization oracle, ensuring it stays within the feasible region \cite{braun2023conditional,Freund2016, pena2023affineinvariantconvergencerates}. {In} \cite{Holloway1974}, Holloway presented a Fully Corrective variant of the algorithm (FCFW). Forty years later, \cite{LacosteJulien2013} gives a linear convergence rate of Away-step Frank-Wolfe for minimizing a smooth and strongly convex objective over a polytope, involving the pyramidal width, a geometry-dependent constant. This result extends to FCFW in \cite{lacostejulien2015global}. The condition measure for the polytope is unified and further clarified in \cite{pena2016polytopeconditioninglinearconvergence}.  This linear convergence generalizes to the sum of a strongly convex function and a linear term \cite{beck2015linearly}. 

\textbf{Dual Correspondance between Frank-Wolfe and Kelley}.
In \cite{bach2013duality}, Bach shows that a generalized version of Frank-Wolfe has a dual correspondence with the mirror descent algorithm. In \cite{Bach2013submodular}, Chapter 7 introduces methods to minimize the sum $f+g$ where $f$ is piecewise linear positive homogeneous. It shows that Kelley is {none other} than FCFW applied to the dual of $g$ on the convex hull of the slopes of $f$. In the same setup, \cite{zhou2018limited} leverages this correspondence and the results from \cite{lacostejulien2015global} to prove the linear convergence of Kelley when $g$ is smooth and strongly convex.

\subsection{Contributions}
We can summarize our contributions as follows:
\begin{enumerate}
    \item We first extend the linear convergence of Kelley's method to a broader class of problems. In particular, we prove that Kelley's method converges linearly for the problem \eqref{problem: primal} with strongly convex and Lipschitz smooth $g$ and general piecewise linear convex $f$, thus extending the previous analysis that requires the additional assumption that $f$ is \textit{positive homogeneous} (see \Cref{section: Linear convergence of Kelley's algorithm}).

    \item We {establish a new duality correspondence between MPB-FA and a Frank-Wolfe (FW) algorithm. This new bridge  allows us to }derive an $\cO(\epsilon^{-1}\log\epsilon^{-1})$ convergence rate of MPB-FA for \eqref{problem: primal} with a convex smooth $g$ and pieceiwse linear $f$, thus improving the best known convergence rate of $\cO(\epsilon^{-2}\log\epsilon^{-1})$ (see \Cref{subsection: Duality of BP}-\Cref{subsection: First convergence rate analysis}). {We believe this new connection will also benefit the study of FW algorithms from the rich literature of the proximal bundle methods.}  

    \item We show through a novel analysis that we can further enhance the convergence rate of the MPB-FA to {$\cO(\epsilon^{-8/9}(\log\epsilon^{-1})^{2/9})$} in the same setting (see \Cref{subsection: Improved convergence rate}). The key arguments of the proof differ significantly from the work on proximal bundle type algorithms \cite{Du2017, Diaz2023optimal,liang2021proximal,LiangMonteiro2023}. The analysis draws upon our dual interpretation of MPB-FA and employs novel geometric arguments. Moreover, it sheds new light on cut management for PBM.  Specifically, it shows the theoretical ground for retaining cuts in the model after a \textit{serious step}. {We also show empirical behavior of MPB-FA and the classical PBM supports the key insights of the improved complexity analysis.} (see \Cref{sec: bundle management}). 
\end{enumerate}

\subsection{Notation and background}
\label{subsec: notations}
We define $\|\cdot\|$ the Euclidean norm and $\langle \cdot , \cdot \rangle$ the corresponding inner product. For a set $\cV$, $\conv(\cV)$ denotes the convex hull of $\cV$. {Let \(\proj_w(\conv(\cV))\) and \(\proj_\beta(\conv(\cV))\) be the projection of $\conv(\cV)$ to the $w$- and $\beta$-space.}
{The convex function $g:\R^n\rightarrow\R$ is $L_g$-smooth on $\R^n$ if the gradient $\nabla g$ satisfies $\|\nabla g(y) - \nabla g(x)\| \leq L_g \|x - y\|$ for all $x, y \in \R^n$. Moreover, $g$ is called $\mu_g$-strongly convex if  $g(y) \geq g(x) + \langle \nabla g(x), y-x \rangle + \frac{\mu_g}{2}\|x - y\|^2$ for all $x, y \in \R^n$.}

For two proper functions $u$, $v$  from $\R^n$ to $\R \cup \{+\infty\}$, let $\square$ denote the infimal convolution over $\R^n$ defined by $(u \,\square\, v)(x)=\inf_{z \in \R^n} \left\{u(x - z) + v(z)\right\}$. 
The Moreau envelope {of $\varphi$ is defined as $\cM_{\rho, \varphi}=\varphi\, \square\, \frac{1}{2\rho}\|\cdot\|^2$, where $\rho > 0$ is the smoothness parameter.} When $\varphi$ is proper, closed, and convex,  $\cM_{\rho, \varphi}$ is a $\frac{1}{\rho}$-smooth convex function that approximates $\varphi$ from below (see Definition 1.22 and Theorem 2.26 in \cite{rockafellar_variational_1998}). 

{For the convex piecewise linear function \(
    f(x) = \max_{(v,b) \in \cV}\{v^\top x + b\}\), let \( D \) denote the \textit{diameter} of \( \cV \), the maximum distance between two points in $\cV$. Let \( \gamma \) be the \textit{pyramidal width} of \(\cV\), defined below, and \( M_f \) the Lipschitz constant of \( f \), i.e., $|f(x)-f(y)|\le M_f\|x-y\|$ for all $x,y \in \R^n$. We further define \(h_k := g+f_k\).}
\leavevmode{\subsubsection{Pyramidal Width}
\label{subsec: pyramidal width}
We use the concept of \emph{pyramidal width} as introduced in \cite{LacosteJulien2013, lacostejulien2015global}. Given a direction $x \in \mathbb{R}^{n}$, the \textit{directional width} of a set $\cV \subseteq \mathbb{R}^{n}$ with respect to $x$ is $\mathop{\mathrm{dirW}}(\cV,\ x) := \max_{v,\ w \mathop{\in} \cV} {(v - w)}^{\top}\frac{x}{\|x\|_2}$.
The \textit{pyramidal width} of $\cV$ is defined as $\mathop{\mathrm{PWidth}}(\cV) \stackrel{\text{def}}{=} \min_{\mathcal{K} \mathop{\in} \mathop{\mathrm{face}(\conv(\cV))}}\min_{x \mathop{\in} \mathop{\mathrm{cone}}(\mathcal{K}-w) \setminus \{0\},\ w \in {\mathcal{K}}}\mathop{\mathrm{PdirW}}(\mathcal{K}\cap \mathop{\cV},\ x,\ w)$,
where $\mathop{\mathrm{PdirW}}$ is the pyramid directional width (see \cite{lacostejulien2015global} for full definition).}
 \vspace{-2mm}
{\subsubsection{Minimum eigenvalue bound}\label{subsec:eigen} In the improved complexity analysis of MPB-FA, we use a result from \cite{kaur_new_2023} on the minimum eigenvalue of a sum of orthogonal projection matrices $aP_1+bP_2$ with positive $a,b$: when $\rank(P_1)+\rank(P_2)=n$ and the first principal angle $\theta_1$ between their ranges is positive, the smallest eigenvalue of $aP_1+bP_2$ is at least $\frac{1}{2}(a+b-\sqrt{(a+b)^2-4ab\sin^2\theta_1})$ (see Theorem 3.5 in \cite{kaur_new_2023}).}

\section{Linear convergence of Kelley's method}
\label{section: Linear convergence of Kelley's algorithm}

This section studies the oracle complexity of the convergence in function value of \Cref{alg: Kelley's algorithm}. This analysis is an extension of \cite{zhou2018limited} in which the function $f$ is \textit{positive homogeneous} (i.e., its epigraph is a polyhedral cone) to any piecewise linear convex function. {While this extension primarily serves as a foundation for \Cref{section : Proximal Bundle algorithm method rate}, we would like to note that some subtlety arises as the dual objective $-g^*(-w)+\beta$ is not strongly convex in $\beta$. In order to apply the existing convergence result on FCFW, we need to estimate a generalized strongly convexity parameter on the dual objective in \Cref{sec:complexity_Kelley}.} We first show in \Cref{subsec: equivalence Kelley FCFW} that \Cref{alg: Kelley's algorithm} and \Cref{alg:FCFW} are dual of one another. \Cref{theorem: linear convergence of Kelley} is the main result of this section and will be used to provide convergence rates for the MPB-FA in \Cref{section : Proximal Bundle algorithm method rate}.

\noindent
\begin{minipage}[t]{0.43\linewidth}
\begin{algorithm}[H]
  \caption{Kelley's Method}\vspace{0.3mm}
  \label{alg: Kelley's algorithm}
\begin{algorithmic}
    \Require \\
    {$f_0(x)$ defined in \eqref{eq: fk} with ${\cV}^{0} \subseteq \cV$.}
    \vspace{2.5mm}
    \For{$k \geq 0$}
      \State Solve subproblem:
      \begin{align}
          x_{k+1}&\in\arg\min_{x \in \R^n} g(x) + f_k(x). \label{eq:Kelley minimization step}
      \end{align}
      \State Add a new cut:
      \vspace{0.2mm}
      \begin{align}
          ({v}_{k+1}, {b}_{k+1}) \in \underset{(v, b) \in \cV}{\argmax}\,\, v^\top  x_{k+1} + b. \label{Kelley's algorithm FW step}
      \end{align}
      \vspace{-2mm}
      \State $\cV^{k+1} \gets \cV^{k}\cup \{({v}_{k+1}, {b}_{k+1})\}.$
      \If{$f(x_{k+1}) - f_k(x_{k+1}) \leq \epsilon$}
        \State Terminate and return $x_{k+1}$. 
      \EndIf
    \EndFor
\end{algorithmic}
\end{algorithm}
\end{minipage}
\hfill
\begin{minipage}[t]{0.55\linewidth}
\begin{algorithm}[H]
  \caption{Fully Corrective FW (FCFW)}
  \label{alg:FCFW}
\begin{algorithmic}
    \Require \\$\widehat{\cV}^{0} \subseteq \cV$, $\phi = g^*(-\cdot)$ a smooth convex function.
    \vspace{2mm}
    \For{$k \geq 0$}
      \State Fully-corrective step:
      \begin{align}
          (w_{k+1}, \beta_{k+1})&\in  \argmin_{(w, \beta) \in \conv(\widehat{\cV}^{k})} \phi(w) - \beta. \label{eq:FCFW fully corrective step}
      \end{align}
      \vspace{-3.5mm}
      \State Frank-Wolfe step:
      \begin{align}
          (v^{FW}_{k+1}, b^{FW}_{k+1}) &\in \underset{(v, b) \in \cV}{\arg\max} -v^\top  \nabla \phi(w_{k+1}) + b.
        \label{FCFW FW step}
      \end{align}
      \vspace{-2.4mm}
      \State $\widehat{\cV}^{k+1} \gets \widehat\cV^{k}\cup \{(v^{FW}_{k+1}, b^{FW}_{k+1})\}.$
          \If{$\langle \nabla \phi(w_{k{+}1}), w_{k{+}1} - v_{k{+}1}^{FW} \rangle - \beta_{k{+}1} + b^{FW}_{k{+}1}  \leq \epsilon$}
        \State Terminate and return $-\nabla\phi(w_{k+1})$. 
      \EndIf
    \EndFor
\end{algorithmic}
\end{algorithm}
\end{minipage}

{\begin{assumption}\label{assum:g strong conv}
Function $g$ is $\mu_g$-strongly convex and $L_g$-smooth.
\end{assumption}}

{In Kelley's method, \Cref{alg: Kelley's algorithm}, at each iteration $k$, a cutting plane model $f_k$ of $f$ is constructed as 
    \begin{align}\label{eq: fk}
        f_k(x) :=\max_{(v,b) \in \cV^{k}}\{v^\top x + b\},
    \end{align} 
where $\cV^{k}\subseteq\cV$ and $f$ is defined in \eqref{eq: f}. Clearly, $f_k$ is a lower approximation of $f$.}

\subsection{Equivalence between \texorpdfstring{\cref{alg: Kelley's algorithm}}{Kelley's method} and \texorpdfstring{\cref{alg:FCFW}}{FCFW}}
\label{subsec: equivalence Kelley FCFW}
{We take three steps to show the equivalence between Kelley's method \cref{alg: Kelley's algorithm} and FCFW \cref{alg:FCFW}. Namely, we prove \eqref{eq:Kelley minimization step} $\iff$ \eqref{eq:FCFW fully corrective step}; \eqref{Kelley's algorithm FW step} $\iff$ \eqref{FCFW FW step}; and finally, the equivalence of the termination conditions of the two algorithms.}
{\subsubsection{Strong duality between \texorpdfstring{\eqref{eq:Kelley minimization step}}{Kelley step 1} and \texorpdfstring{\eqref{eq:FCFW fully corrective step}}{FCFW step 1}} \;}

\begin{lemma}[Lagrangian dual {and strong duality}]
\label{lm: Kelley's dual problem}
The dual of \eqref{problem: primal} is
\begin{align}
        D(\cV):\quad \underset{(w, \beta) \in \conv(\cV) }{\max} \, -g^*(-w) + \beta =- \underset{(w, \beta) \in \conv(\cV) }{\min} \, \phi(w) - \beta,\label{problem: Dual}
\end{align}
where $\phi(w) := g^*(-w)$.  {Strong duality holds between \eqref{problem: primal} and \eqref{problem: Dual}. When $\cV^k=\widehat{\cV}^k$ in iteration $k$ of the two algorithms, strong duality holds between \eqref{eq:FCFW fully corrective step} and \eqref{eq:Kelley minimization step}.}
\end{lemma}
\begin{proof}
{The primal problem \eqref{problem: primal} can be written as 
\begin{align}
     \min_{x\in\R^n} \; g(x)+f(x) & = \min_{x, z}\biggl\{g(x) + z \; : \; z\ge v^\top x + b, \;\; \forall (v,b)\in \cV\biggr\}.\label{eq:primal_pwlconstr}
\end{align}
Introducing multipliers $\lambda_{v,b}$ for constraints in \eqref{eq:primal_pwlconstr}, the Lagrangian dual is}
{ 
\begin{align}
    & \max_{\lambda \ge 0} \; \min_{x,z} \; g(x)+z\biggl(1-\sum_{(v,b) \in \cV}\lambda_{v,b}\biggr) + \biggl(\sum_{(v,b) \in \cV}\lambda_{v,b}  v\biggr)^\top x  + \sum_{(v,b)\in\cV} \lambda_{v,b} b \notag\\
   = & \max_{\lambda \ge 0, e^\top \lambda = 1} \; \min_{x} \; g(x) + \biggl(\underbrace{\sum_{(v,b) \in \cV}\lambda_{v,b}  v}_{=: w}\biggr)^\top x  + \underbrace{\sum_{(v,b)\in\cV} \lambda_{v,b} b}_{=:\beta}\label{eq:dual_min_z}\\
    =& \max_{(w,\beta)\in\conv(\cV)}\min_{x}\; g(x) + w^\top x + \beta = \max_{(w,\beta)\in\conv(\cV)} -g^*(-w)+\beta,\label{eq:dual_w,beta}
\end{align}
where \eqref{eq:dual_min_z} is obtained by minimizing over $z$ first and $e$ is the vector with all $1$'s; then $\lambda$ represents the weights of a convex combination. Then, $(w,\beta)$ are in the convex hull of $\cV$. Lastly, we use the conjugate of $g$. If the sets of cuts in the two algorithms are the same, i.e., $\cV^k=\widehat{\cV}^k$, then by the same argument, the dual of \eqref{eq:Kelley minimization step} is \eqref{eq:FCFW fully corrective step}}. 
{In \eqref{eq:dual_w,beta}, $g(x)+w^\top x + \beta$ is convex in $x$ and linear in $(w,\beta)$. Furthermore, \(\conv(\cV)\) is a nonempty convex compact set due to the piecewise linearity of $f$. The Sion's minimax theorem \cite{Sion1958} applies and strong duality holds between  \eqref{eq:dual_w,beta} and \eqref{problem: primal}. Therefore, strong duality holds between \eqref{problem: primal} and \eqref{problem: Dual} and between \eqref{eq:Kelley minimization step} and \eqref{eq:FCFW fully corrective step}.}
\end{proof}
{Immediately, we have the following corollary.
\begin{corollary}[Unique correspondence between $x_{k+1}$ and $(w_{k+1},\beta_{k+1})$]\label{cor:kelley sol correspond}
When Algorithms \ref{alg: Kelley's algorithm} and \ref{alg:FCFW} have the same cuts $\cV^k=\widehat{\cV}^k$, \eqref{eq:Kelley minimization step} and \eqref{eq:FCFW fully corrective step} have a unique pair of optimal solutions $x_{k+1}$ and $(w_{k+1},\beta_{k+1})$.
\end{corollary}}
\subsubsection{{Equivalence between \texorpdfstring{\eqref{Kelley's algorithm FW step}}{Kelley2} and \texorpdfstring{\eqref{FCFW FW step}}{FCFW2}}} {This step goes deeper to describe a precise relation} between $x_k$ in the primal \Cref{alg: Kelley's algorithm} and $(w_k, \beta_k)$ in the dual \Cref{alg:FCFW}, which will imply the equivalence between  \texorpdfstring{\eqref{Kelley's algorithm FW step}}{Kelley2} and \texorpdfstring{\eqref{FCFW FW step}}{FCFW2}.

\begin{lemma}[{Explicit relation} between primal and dual solutions]
\label{lemma:Primal-Dual variables correspondence}
{Assume \Cref{assum:g strong conv}.} 
Let $(w_{k+1}, \beta_{k+1})$ be optimal for the dual minimization step \eqref{eq:FCFW fully corrective step} and $x_{k+1}$ be the primal optimal solution of \eqref{eq:Kelley minimization step}. {If $\cV^k=\widehat{\cV}^k$}, then
\begin{align}
    w_{k+1} &= -\nabla g(x_{k+1}), \label{eq:primal_dual_w}\\
     x_{k+1} &= \nabla g^*(-w_{k+1}).\label{eq:primal_dual_x}
\end{align}
\end{lemma}
\begin{proof}
\leavevmode{We prove in two steps.
Step 1: Lemma \ref{lm: Kelley's dual problem} establishes strong duality
\begin{align}
    \min_{x\in\R^n}\biggl\{g(x) + \underbrace{\max_{(w,\beta)\in\conv(\cV^k)} w^\top x + \beta}_{f_k(x)}\biggr\} = \max_{(w,\beta)\in\conv(\cV^k)}\min_{x\in\R^n} g(x) + w^\top x + \beta,\label{eq:proof:strong_dual}
\end{align}
where the left hand side (LHS) is \eqref{eq:Kelley minimization step} and the right hand side (RHS) is \eqref{eq:FCFW fully corrective step} when $\cV^k=\widehat\cV^k$. Hence, $x_{k+1}$ is optimal for the outer min on LHS of \eqref{eq:proof:strong_dual} and $(w_{k+1}, \beta_{k+1})$ is optimal for the outer max on RHS. By Saddle Point Theorem \cite{BorweinLewis}, both sides of \eqref{eq:proof:strong_dual} are equal to $g(x_{k+1})+w_{k+1}^\top x_{k+1}+\beta_{k+1}$. 
Therefore, $x_{k+1}$ is optimal for the inner min of the RHS of \eqref{eq:proof:strong_dual}. Thus, $\nabla g(x_{k+1})+w_{k+1}=0$, proving \eqref{eq:primal_dual_w}. 

Step 2: Using conjugacy, $\min_{x\in\R^n}\{g(x)+w_{k+1}^\top x + \beta_{k+1}\}=-g^*(-w_{k+1})+\beta_{k+1}$. Since $g$ is strongly convex, $\nabla g^*$ exists and the inner min on the RHS of \eqref{eq:proof:strong_dual} has a unique optimum $x_{k+1}$ by Step 1. Thus, $\nabla (-g^*(-w_{k+1}))=x_{k+1}$, proving \eqref{eq:primal_dual_x}.}
\end{proof}
\leavevmode{
This immediately gives the following corollary.
\begin{corollary}\label{corr:kelley=FCFW}
If $\cV^k=\widehat\cV^k$, then \texorpdfstring{\eqref{Kelley's algorithm FW step}}{Kelley2} and \texorpdfstring{\eqref{FCFW FW step}}{FCFW2} are equivalent. That is,
    \begin{align}
        \arg\max_{(v,b)\in\cV} v^\top x_{k+1}+b = \arg\max_{(v,b)\in\cV} -v^\top \nabla \phi(w_{k+1})+b. \label{eq:cor:KelleyFCFW}
    \end{align}
\end{corollary}
}
{\subsubsection{Equivalence of the termination conditions} First, we need the following definition.}
\begin{definition}[Primal gap and Frank–Wolfe gap in Definitions 1.10-11 in \cite{braun2023conditional}]
\label{def: FW gap}
Let $s : \conv(\cV) \to \R^n$ {be convex and smooth} and $x^* \in \arg\min_{x \in \conv(\cV)}\, s(x)$. The primal gap is defined as $d^{primal}(x) := s(x) - s(x^*)$.  The Frank–Wolfe gap (aka the dual gap) of $s$ at $x$ is $
d^{FW}(x) :=  \max_{v\in\cV} \langle\nabla s (x), x - v\rangle$. 
\end{definition}
{\begin{remark}\label{rm:primalgap_FWgap}
Clearly, $d^{primal}(x)\ge 0$. By convexity of $s$, we also have $d^{FW}(x) \geq 0$:
\begin{align}
 d^{FW}(x) = \max_{v\in\cV} \;\langle\nabla s (x), x - v\rangle \ge  \langle\nabla s (x), x - x^*\rangle \ge d^{primal}(x) \ge 0. \label{ineq: primal gap < FW gap}
\end{align}
\end{remark}}

The following lemma shows the correspondence between the termination criterion of Algorithms \ref{alg: Kelley's algorithm} and \ref{alg:FCFW}. 

\begin{lemma}
\label{lemma: test in Kelley corresponds to Frank Wolfe Gap} 
{If $\cV^k=\widehat\cV^k$}, the termination criterion in \Cref{alg: Kelley's algorithm} corresponds to the termination criterion in \Cref{alg:FCFW}, which   {is} a comparison between the accuracy $\epsilon$ and the Frank-Wolfe gap associated with the dual problem \eqref{problem: Dual}. 
\end{lemma}
\begin{proof}
 {If $\cV^k=\widehat\cV^k$, then}
\begin{align}
    f(x_{k+1}) - f_k(x_{k+1}) &:=
     \underset{(v,b) \in \cV}{\max} \,\left\{v^\top  x_{k+1} + b\right\} - \underset{(w,\beta) \in \conv(\cV^{k})}{\max} \,\left\{w^\top  x_{k+1} + \beta\right\} \nonumber\\ 
     &= (v^{FW}_{k+1})^\top  x_{k+1} +b^{FW}_{k+1} -\left( w_{k+1}^\top  x_{k+1} +\beta_{k+1}\right) \label{maximum on Vt replaced by wt+1} \\
    &=\langle \nabla \phi(w_{k+1}), w_{k+1} - v_{k+1}^{FW} \rangle - \beta_{k+1} + b^{FW}_{k+1} \label{eq: recover gap from FCFW}\\
   &= \biggl\langle  \binom{\nabla\phi(w_{k+1})}{-1}, \binom{w_{k+1}}{\beta_{k+1}} - \binom{v_{k+1}^{FW}}{b^{FW}_{k+1}}\biggr\rangle.
   \label{kelley Null step condition is FW Dual gap} 
\end{align}
We get \eqref{maximum on Vt replaced by wt+1} from \eqref{eq:FCFW fully corrective step} and \eqref{FCFW FW step}. In \eqref{eq: recover gap from FCFW} we replace $x_{k+1}$ by $-\nabla \phi(w_{k+1})$ using \Cref{lemma:Primal-Dual variables correspondence} and remark that it is exactly the termination criterion of \Cref{alg:FCFW}.  {Equation} \eqref{kelley Null step condition is FW Dual gap} is the Frank-Wolfe gap (\Cref{def: FW gap}) associated with problem \eqref{problem: Dual} at $(w_{k+1}, \beta_{k+1})$.  
\end{proof}

 {Finally, we} establish that  {if Algorithms \ref{alg: Kelley's algorithm} and \ref{alg:FCFW} start from the same initial cuts $\cV^0=\widehat\cV^0$}, then all their iterations have an exact correspondence, which allows us to directly translate the convergence results of \Cref{alg:FCFW} to \ref{alg: Kelley's algorithm}.
\begin{proposition}
    \label{lm:correspondence Kelley FCFW}
     {Under \Cref{assum:g strong conv} and $f$ given in \eqref{eq: f}, when $\cV^0=\widehat\cV^0$,} the iterates of \Cref{alg: Kelley's algorithm} match the iterates of \Cref{alg:FCFW}  {in the sense that for each iteration $k$, 1) $x_{k+1}$ and $(w_{k+1},\beta_{k+1})$ are the unique optimal solutions to the primal-dual pair \eqref{eq:Kelley minimization step} and \eqref{eq:FCFW fully corrective step}; 2) problems \texorpdfstring{\eqref{Kelley's algorithm FW step}}{Kelley2} and \texorpdfstring{\eqref{FCFW FW step}}{FCFW2} are exactly the same, so their optimal solutions can be chosen so that $(v_{k+1},b_{k+1})=(v^{FW}_{k+1}, b^{FW}_{k+1})$, which guarantees $\cV^{k+1}=\widehat\cV^{k+1}$; and 3) the two algorithms terminate at the same iterate}.     
\end{proposition}
\begin{proof} 
 {Since both algorithms start from the same cuts $\cV^0=\widehat\cV^0$ and \Cref{assum:g strong conv} holds, by \Cref{lemma:Primal-Dual variables correspondence}, $x_1$ and $(w_1, \beta_1)$ are the unique solutions to \eqref{eq:Kelley minimization step} and \eqref{eq:FCFW fully corrective step}, 
and \eqref{Kelley's algorithm FW step} is \eqref{FCFW FW step} by \Cref{corr:kelley=FCFW}. Thus,  \(({v}_{1}, {b}_{1})=(v^{FW}_{1}, b^{FW}_{1}) \implies \cV^1=\widehat\cV^1\).  Recursively, $\cV^k=\widehat\cV^k$ for $k\ge 0$. By \Cref{lemma: test in Kelley corresponds to Frank Wolfe Gap}, the two algorithms  terminate at the same iterate.}
\end{proof}


\subsection{Iteration complexity of \texorpdfstring{\Cref{alg: Kelley's algorithm}}{FCFW}}\label{sec:complexity_Kelley}
In the setup of \cite{zhou2018limited}, where the dual problem is smooth and strongly convex, FCFW on a polytope converges linearly \cite{lacostejulien2015global}.  {We generalize the results in \cite{zhou2018limited} to problem  \eqref{problem: Dual}, where the objective function $-g^*(-w)+\beta$ is not strongly convex in \((w,\beta)\) because $g^*(-w)$ is not a function of $\beta$. In general, without strong convexity, Kelley or FCFW does not enjoy linear rate. However, for special structured problems such as the sum of a linear function and a strongly convex function composed with a linear map over a polytope, the literature  has derived linear convergence rate for Frank-Wolfe type algorithms \cite{lacostejulien2015global, beck2015linearly}. } 
The proof of the linear convergence of FCFW in \cite{lacostejulien2015global} directly uses a convergence result established for the Away-Step Frank-Wolfe (AFW) in \cite{beck2015linearly}. In \Cref{lemma: strong convexity constant of psi}, we provide a geometric strong convexity constant $\bar{\mu}_\psi$ for \eqref{problem: Dual}. Compared to \cite{beck2015linearly}, we eliminate the dependence on the Hoffman constant that may be arbitrarily large. As in \cite{lacostejulien2015global}, we will extend the obtained convergence rate for AFW to \Cref{alg:FCFW} by noticing that a step of FCFW makes at least as much progress as the corresponding descent step and all the following away steps of AFW. This will lead to \Cref{theorem: linear convergence of Kelley} for linear convergence of \Cref{alg: Kelley's algorithm}.


By the exact correspondence between Kelley and FCFW in \Cref{lm:correspondence Kelley FCFW}, to obtain convergence rate for Kelley on \eqref{problem: primal}, it suffices to analyze the convergence rate of FCFW applied to \eqref{problem: Dual}, which is written again for convenience
\begin{align}
    \min_{(w,b) \in  {\conv(\cV)}} \,  {\psi(w,\beta)\coloneqq g^*(-w)-\beta}. \label{eq:dual_restated}
\end{align}

 {We first obtain a generalized geometric strong convexity constant of $\psi$. Leveraging this constant, we can directly apply the result in \cite{lacostejulien2015global} to our problem.}  {From now on, we assume $\cV^k=\widehat\cV^k$ for all $k$ and use $\cV^k$ in place of $\widehat\cV^k$.}

\subsubsection{ {Generalized geometric strong convexity of the dual objective}}\label{subsec:generalized_strong_convexity}

\begin{lemma}
\label{lemma: strong convexity constant of psi}
 {Under \Cref{assum:g strong conv} and $f$ being $M_f$-Lipschitz continuous, let $D_b$ be the diameter of \(\proj_\beta(\conv(\cV))\).} Let $(w^*,\beta^*)$ be the unique optimal solution of  \eqref{eq:dual_restated}.
 {Then the dual objective $\psi(w,\beta)$ defined in \eqref{eq:dual_restated}}  satisfies for all  $(w,\beta)\in \cV $:
    \begin{align*}
         \psi(w^*,\beta^*) - \psi(w,\beta) - 2 \langle \nabla \psi (w, \beta), ((w^*)^\top , \beta^*)^\top  {-} (w^\top , \beta)^\top  \rangle \geq 2\bar{\mu}_\psi \left\|\binom{w^*}{ \beta^*} - \binom{w}{ \beta} \right\|^2
    \end{align*}
where the generalized geometric strong convexity constant \(\bar{\mu}_\psi\) of $\psi$ is given by
    \begin{align}
        \bar{\mu}_\psi =\frac{1}{2}\left( D_b + 3 \frac{ {4} M_f^2}{\mu_g}  +6M_f \|x^*\| + 2L_g \left[\left(\frac{ {2}M_f}{\mu_g} + \|x^*\|\right)^2 + 1\right]\right)^{-1}.\label{eq: generalized strong convexity constant}
    \end{align}
\end{lemma}
The proof can be found in the Appendix.


\subsubsection{ {Linear convergence of Kelley's method for the primal}}\label{subsec:linear_conv_Kelley} 

\begin{theorem}[Linear convergence of Kelley's method]
    \label{theorem: linear convergence of Kelley}
    Under \Cref{assum:g strong conv}, 
    \Cref{alg: Kelley's algorithm} applied to \eqref{problem: primal} finds an $\epsilon$-optimal solution in at most
    \begin{align}
        1 + \max\left\{2, \frac{D^2}{\bar{\mu}_\psi\mu_g \gamma^2}\right\} \log \left(\frac{D^2}{2 \epsilon \mu_g} \right)
    \end{align}
    iterations, where $\bar{\mu}_\psi$ is given by  {\eqref{eq: generalized strong convexity constant}, $D$ is the diameter of $\conv(\cV)$}, and $\gamma$ is the pyramidal width of $\cV$.
\end{theorem}
\begin{proof}
    We leverage the correspondence between Kelley's method and FCFW shown in \Cref{lemma:Primal-Dual variables correspondence}, thereby bounding the iteration complexity of Kelley's method. Indeed, Theorem 8 in \cite{lacostejulien2015global} has shown the upper bound of FCFW oracle complexity stated in \Cref{theorem: linear convergence of Kelley} above, with the constant $\bar{\mu}_\psi$ as defined in \Cref{lemma: strong convexity constant of psi}, a lower bound of the generalized geometric strong convexity constant of $\psi$.  
\end{proof}


\section{Modified Proximal Bundle with Fixed Absolute Accuracy}
\label{section : Proximal Bundle algorithm method rate}
 {This section studies problem \eqref{problem: primal} \(\min_{x\in\R^n} \; g(x)+f(x)\)
under the following more general assumption.} 
 {\begin{assumption}\label{assum: general}
Function $g:\R^n\to\R$ is convex, $L$-smooth, but not necessarily strongly convex. The function $f$ is given by \eqref{eq: f}.
\end{assumption}}

  In this setup, we turn our attention to the proximal bundle methods that add strong convexity to $g$ through a quadratic proximal term. We   {propose the modified proximal bundle method with fixed absolute error (MPB-FA) in \Cref{alg: Proximal Bundle Method}. The dual \Cref{alg:dual proximal bundle} turns out to be a smoothed FCFW, dealing with the nonsmooth nature of $g^*$.} We study the convergence rate of MPB-FA by taking advantage of this duality correspondence.  \Cref{subsection: Duality of BP} proves the correspondence between \Cref{alg: Proximal Bundle Method} and \Cref{alg:dual proximal bundle}. 
   {\Cref{subsec: Comparison between MPB-FA and PBM} points out the main difference between MPB-FA and the classical PBM in terms of the serious step test and discusses the interesting affine invariance property of MPB-FA.  \Cref{subsection: First convergence rate analysis} gives a first convergence analysis of MPB-FA by leveraging the duality correspondence and \Cref{theorem: linear convergence of Kelley}.  \Cref{subsection: Improved convergence rate} presents a novel and deep analysis of MPB-FA and obtains} 
 \Cref{thm: 1/e convergence when g is smooth}, which achieves an  {$\cO(\epsilon^{-8/9})$} convergence rate, improving the rate in \Cref{subsection: First convergence rate analysis}.

\begin{figure}[h!]
\centering
\noindent
\begin{minipage}[t]{0.46\linewidth}
\begin{algorithm}[H]
  \caption{MPB-FA}
  \label{alg: Proximal Bundle Method}
\begin{algorithmic}
\Require \\
$x_0 \in \R^n$, $\cV^{0} \subset \cV$, $\epsilon$, $\rho>0$
    \For{$k \geq 0$}
    \State Compute $y_{k + 1}$ by solving :
      \begin{align}
         \min_{x \in \R^n} g(x) {+} f_k(x) {+} \frac{\rho}{2}\|x {-} x_k\|^2 \label{eq:proxbundle_alg}
      \end{align}
    \vspace{-5mm}
      \State \begin{align}({v}_{k+1}, {b}_{k+1}) \in \arg\max_{(v, b) \in \cV}\; (v^{T} y_{k+1} + b)\label{bundle FW step}
      \end{align}
       \vspace{-2mm}
      \State Update $\cV^{k+1} \gets \cV^{k} \cup \{({v}_{k+1}, {b}_{k+1})\}$ 
      \State   {Serious step test} : 
      \If{$f(y_{k+1}) - f_k(y_{k+1}) \leq \epsilon/2$} 
        \State $x_{k+1} \gets y_{k+1}$ \Comment{(serious step)}
      \Else
        \State $x_{k+1} \gets x_k$ \Comment{(null step)}
      \EndIf
    \EndFor
\end{algorithmic}
\end{algorithm}
\end{minipage}
\hfill
\begin{minipage}[t]{0.52\linewidth}
\begin{algorithm}[H]
  \caption{Modified FCFW}
  \label{alg:dual proximal bundle}
  \begin{algorithmic}
    \Require \, $\widehat{x}_0 \in \R^n$, $\widehat\cV^{0} \subset \cV$, $\epsilon/2$, $\rho>0$
    \For{$k \geq 0$}
      \State Compute $({w}_{k+1}, {\beta}_{k+1})$ by solving :
      \vspace{-1.5mm}
        \begin{align}
                \underset{(w, \beta) \in \conv(\widehat\cV^k)}{\min} \, \cM_{\rho, \phi} (w - \rho \widehat{x}_k) - \beta   \label{dual bundle method fully corrective step} 
        \end{align}
        \vspace{-5mm}
        \begin{align}
            \widehat{y}_{k+1} \gets -\nabla \cM_{\rho, \phi}({w}_{k+1} - \rho \widehat{x}_k) \label{dual bundle method first order oracle}
        \end{align}
        \vspace{-5.5mm}
    \begin{align}
        (v_{k+1}^{FW}, b_{k+1}^{FW}) \in\arg\max_{(v, b) \in \cV}\; (v^{T} \widehat{y}_{k+1} + b)\label{eq: dual bundle method FW step}
    \end{align}
      \State Update $\widehat\cV^{k+1} \gets \widehat\cV^{k} \cup \{(v_{k+1}^{FW}, b_{k+1}^{FW})\}$ 
      \vspace{3mm}
      \If{$\langle \widehat{y}_{k+1}, v_{k+1}^{FW}  - w_{k+1} \rangle + b^{FW}_{k+1} -\beta_{k+1}  \leq \epsilon/2$}
           \State $\widehat{x}_{k+1} \gets \widehat{y}_{k+1}$
      \Else
        \State $\widehat{x}_{k+1} \gets \widehat{x}_k$
      \EndIf 
    \EndFor
    \end{algorithmic}
\end{algorithm}
\end{minipage}
\raggedright
{  {Recall that \(f_k(x)=\max_{(v,b)\in\cV^k} v^\top x+b\), and $\cM_{\rho, \phi}(w):=\min_v \{\phi(v) + \frac{1}{2\rho}\|v-w\|^2\}$ stands for the Moreau envelope of the function $\phi\coloneqq g^*(-\cdot)$ with parameter $\rho$ as defined in \Cref{subsec: notations}. 
}}
\end{figure}
 
\subsection{Equivalence between \texorpdfstring{\Cref{alg: Proximal Bundle Method}}{PBM} and \texorpdfstring{\Cref{alg:dual proximal bundle}}{Modified FCFW}}
\label{subsection: Duality of BP}
 {Similar to \Cref{subsec: equivalence Kelley FCFW}, we take three steps to show  \eqref{eq:proxbundle_alg} $\iff$ \eqref{dual bundle method fully corrective step}; \eqref{bundle FW step} $\iff$ \eqref{eq: dual bundle method FW step}; and the serious step test in \Cref{alg: Proximal Bundle Method} $\iff$ the corresponding ``if'' test in \Cref{alg:dual proximal bundle}. }

 {\subsubsection{Strong duality between \texorpdfstring{\eqref{eq:proxbundle_alg}}{bundle} and \texorpdfstring{\eqref{dual bundle method fully corrective step}}{dual bundle}} The proof largely follows that of \Cref{lm: Kelley's dual problem} with the twist of handling regularized version of $g$.}

\begin{lemma}
\label{lm: dual of prox bundle minimization}
Suppose  {$\cV^k=\widehat\cV^k$} at iteration $k$. 
The dual of problem \eqref{eq:proxbundle_alg} is 
\begin{align}
    \tilde{\mathcal{D}}(\cV^k, x_k) : 
    -\min_{(w, \beta)\in\conv({\cV}^k)} \left\{\cM_{\rho, g^*}(\rho x_k  - w) - \beta\right\} + \frac{\rho}{2}\|x_k\|^2. \label{eq:proximal bundle dual}
\end{align}
\end{lemma}
\begin{proof}
 {Define $\hat{g}_k \coloneq g + \frac{\rho}{2}\|\cdot - x_k\|^2$. Then \eqref{eq:proxbundle_alg} can be written as \(\min_{x\in\R^n} \hat{g}(x)+f_k(x)\), which has the same form as \eqref{eq:Kelley minimization step} with $\hat{g}$ replacing $g$. By \Cref{lm: Kelley's dual problem}, the dual of \eqref{eq:proxbundle_alg} is \(\min_{(w,\beta)\in\cV} \hat{g}^*(-w)-\beta\). All we need to do is to compute the conjugate of $\hat{g}$:}
\begin{align*}
    \hat{g}_k^*(-w) &= \left(g + \frac{\rho}{2}\|\cdot - x_k\|^2\right)^*(-w)= \left(g^* \square \left(\frac{\rho}{2}\|\cdot - x_k\|^2\right)^* \right)(-w)
    \\&= \left(g^* \square (\frac{1}{2\rho}\|\cdot\|^2 + \langle \cdot , x_k\rangle) \right)(-w)
    = \left(g^* \square (\frac{1}{2\rho}\|\cdot + \rho x_k\|^2) \right)(-w) - \frac{\rho}{2}\|x_k\|^2\\
    &=\cM_{\rho, g^*}(\rho x_k  - w) - \frac{\rho}{2}\|x_k\|^2.\stepcounter{equation}\tag{\theequation}\label{eq: gk conjugate is moreau enveloppe}
\end{align*}
Leaving out the constant term $-\rho/2\|x_k\|^2$, we get the dual \eqref{eq:proximal bundle dual}.  {To recover \eqref{dual bundle method fully corrective step}, we note that $\phi(w):=g^*(-w)$. Thus, $\cM_{\rho, g^*}(\rho x_k  - w)=\cM_{\rho, \phi}(w - \rho x_k)$. Moreover, strong duality holds between \eqref{eq:proxbundle_alg} and \eqref{dual bundle method fully corrective step}  due to the compactness of $\conv(\cV^k)$, following a similar argument as in \Cref{lm: Kelley's dual problem}}. 
\end{proof}
 {Immediately, we have the following corollary.
\begin{corollary}\label{cor:bundl_y_w}
When Algorithms \ref{alg: Proximal Bundle Method} and \ref{alg:dual proximal bundle} have the same cuts $\cV^k=\widehat{\cV}^k$, the primal \eqref{eq:proxbundle_alg} and the dual \eqref{dual bundle method fully corrective step} have  unique optimal solution $y_{k+1}$ and $(w_{k+1},\beta_{k+1})$.
\end{corollary}}
\vspace{-0.5cm}
 {\subsubsection{Equivalence between \texorpdfstring{\eqref{bundle FW step}}{bundle} and \texorpdfstring{\eqref{eq: dual bundle method FW step}}{dual bundle}}  
Now we can show the precise relation between $y_{k+1}$ and $(w_{k+1},\beta_{k+1})$. }
\begin{lemma}
\label{lm: dual correspondence variables BP}
 {Under \Cref{assum: general}, if $\cV^k=\widehat\cV^k$ and $x_k=\widehat{x}_k$ in Algorithms \ref{alg: Proximal Bundle Method} and \ref{alg:dual proximal bundle},}
then $y_{k+1}$ and $(w_{k+1}, \beta_{k+1})$ are the unique optimal solution to \eqref{eq:proxbundle_alg} and \eqref{dual bundle method fully corrective step}. 
Moreover,  {$y_{k+1}=\widehat{y}_{k+1}$}. Thus,
\begin{align}
    {y}_{k+1} = -\nabla \cM_{\rho, \phi}(w_{k+1} - \rho x_k). \label{eq: definition of y in terms of grad of moreau enveloppe}
\end{align}
\end{lemma}\label{lm:precise_relation_bundle}
\begin{proof}
 {From \Cref{lm: dual of prox bundle minimization}, we know  the conjugate of $\hat{g}_k \coloneq g + \frac{\rho}{2}\|\cdot - x_k\|^2$ is the Moreau envelope \(\cM_{\rho, g^*}(\rho x_k  - w)\). Since $\hat{g}_k$ is strongly convex and smooth, the Moreau envelope is also strongly convex and smooth. If \eqref{eq:proxbundle_alg} and \eqref{dual bundle method fully corrective step} share the same set of cuts and proximal center, i.e., $\cV^k=\widehat\cV^k$ and $x_k=\widehat{x}_k$, then they are a primal-dual pair, and both have a unique optimal solution.} 
Applying \Cref{lemma:Primal-Dual variables correspondence} by replacing  $g$ by $\hat{g}_k$ and using the above Moreau envelope as the conjugate of $\hat{g}_k$ leads to \eqref{eq: definition of y in terms of grad of moreau enveloppe}. Using the definition \eqref{dual bundle method first order oracle}, we have $y_{k+1}=\widehat{y}_{k+1}$.
\end{proof}
 {Immediately, we have the following corollary.
\begin{corollary}
\label{cor:bundle_step2}
If $\cV^k=\widehat\cV^k$ and $x_k=\widehat{x}_k$ in iteration $k$ of \Cref{alg: Proximal Bundle Method} and \Cref{alg:dual proximal bundle}, then \eqref{bundle FW step} and \eqref{eq: dual bundle method FW step} are the same problem with $y_{k+1}=\widehat{y}_{k+1}$. Thus, $    \arg\max_{(v, b) \in \cV}\; (v^{T} y_{k+1} + b) = \arg\max_{(v, b) \in \cV}\; (v^{T} \widehat{y}_{k+1} + b)$.
\end{corollary}}

\subsubsection{Equivalence between the tests}
For iteration $k$ of \Cref{alg:dual proximal bundle}, we define $\Psi_k(w,\beta) := \cM_{\rho, \phi}(w - \rho x_k) - \beta$, then
$\nabla \Psi_k(w, \beta) = \binom{\nabla \cM_{\rho, \phi}(w - \rho x_k))}{-1}$. Note that problem \eqref{dual bundle method fully corrective step} can be equivalently written as
\begin{align}
    \underset{(w,\beta) \in \conv(\cV^k)}{\min}\, \Psi_k(w,\beta). \label{problem : null step dual problem}
\end{align}
The following \Cref{lm:bundle_tests} shows that the serious step test in \Cref{alg: Proximal Bundle Method} is chosen to correspond to the Frank-Wolfe gap in \Cref{alg:dual proximal bundle}. 
 {\begin{lemma}
\label{lm:bundle_tests}
\label{lemma: tests equivalence} 
 {If $\cV^k=\widehat\cV^k$ and $x_k=\widehat{x}_k$}, the serious step test in \Cref{alg: Proximal Bundle Method} is the same as the ``if'' test in \Cref{alg:dual proximal bundle}, which is the Frank-Wolfe gap of \eqref{dual bundle method fully corrective step}.
\end{lemma}}
\begin{proof}
     {If $\cV^k=\widehat\cV^k$ and $x_k=\widehat{x}_k$, then, by \Cref{lm: dual correspondence variables BP}, $y_{k+1}=\widehat{y}_{k+1}$. Thus, 
\begin{align}
    f(y_{k+1}) - f_k(y_{k+1}) 
     &=\underset{(v,b) \in \cV}{\max} \,\left\{v^\top  \widehat{y}_{k+1} + b\right\} - \underset{(w,\beta) \in \conv(\widehat\cV^{k})}{\max} \,\left\{w^\top  \widehat{y}_{k+1} + \beta\right\} \nonumber\\
     &= (v^{FW}_{k+1})^\top  \widehat{y}_{k+1} +b^{FW}_{k+1} -\left( w_{k+1}^\top  \widehat{y}_{k+1} +\beta_{k+1}\right), \label{proof_bundle_test_FCFW}\\
     &= \left\langle \nabla \Psi_k(w_{k+1}, \beta_{k+1}), \binom{w_{k+1}}{\beta_{k+1}} - \binom{v_{k+1}^{FW}}{b^{FW}_{k+1}}\right\rangle,
   \label{Null step condition is FW Dual gap} 
\end{align}
where \eqref{proof_bundle_test_FCFW} establishes the equivalence of the two tests and \eqref{Null step condition is FW Dual gap} further gives the Frank-Wolfe gap of \eqref{dual bundle method fully corrective step} due to \eqref{eq: dual bundle method FW step} and \Cref{def: FW gap}.}
\end{proof}

  {Finally, we come to the equivalence of \Cref{alg: Proximal Bundle Method} and \Cref{alg:dual proximal bundle}.
\begin{proposition}
\label{prop:equivalence_bundle}
     {Under \cref{assum: general}, when $\cV^0=\widehat\cV^0$ and $x_0=\widehat{x}_0$,} the iterates of \Cref{alg: Proximal Bundle Method}  {can be chosen such that they exactly} match the iterates of \Cref{alg:dual proximal bundle}  and they terminate at the same iterate.
    \end{proposition}}
\begin{proof}  {
    The proof follows by combining \Cref{cor:bundl_y_w}, \Cref{cor:bundle_step2}, and \Cref{lm:bundle_tests}, and using the same arguments as in the proof of \Cref{lm:correspondence Kelley FCFW}.}
\end{proof}

\subsection{Comparison between MPB-FA and PBM}
\label{subsec: Comparison between MPB-FA and PBM}
 {The key difference between MPB-FA, \Cref{alg: Proximal Bundle Method}, and the classical PBM is the serious step test. In particular, MPA-FA has 
the \textit{serious step test}} as $f(y_{k+1}) - f_k(y_{k+1}) \leq \epsilon/2$, while the classical PBM \cite{Lemarechal1995, Du2017, Diaz2023optimal} has a serious step test as 
\begin{align}
    \beta(f(x_{k}) - f_k(y_{k+1})) \leq f(x_{k}) - f(y_{k+1}),\label{eq:classicPBA_test}
\end{align} 
with $\beta \in (0,1)$, 
which compares the progress made between the iterates $x_k$ and $y_{k+1}$ to the expected progress. 

The  {serious step} test in \Cref{alg: Proximal Bundle Method} is inspired by the test introduced in \cite{liang2021proximal,LiangMonteiro2023}. In the classical setup where $g=0$, their test is $f(z_{k+1}) - f_k(y_{k+1}) \leq \frac{\rho}{2}\|y_{k+1}-x_k\|^2+\epsilon/2$, where $z_{k+1}$ represents the best iterate for the proximal problem.  Within the context of the composite problem \eqref{problem: primal}, these tests maintain the algorithm's invariance with respect to the redistribution of a linear map between $f$ and $g$.   {Indeed, for any $a \in \mathbb{R}^n$, we replace $g$ by $\hat{g} : x \mapsto g(x) - a^\top x$ and $f$ by $\hat{f} : x \mapsto f(x) + a^\top x$. The new piecewise linear function $\hat{f}$ is the supremum of the affine functions defined by $\hat{\mathcal{V}} \coloneqq \{ (v + a, b) : (v, b) \in \mathcal{V} \}.$
A simple recursive argument shows that for the modified objective $\hat{f} + \hat{g}$, one has $\hat{v}_k = v_k + a$ and $\hat{b}_k = b_k$, implying that $\hat{f}_k(x) = f_k(x) + a^\top x$ for all $x \in \mathbb{R}^n$. The minimization step in \Cref{eq:proxbundle_alg} remains unchanged since $g + f_k = \hat{g} + \hat{f}_k$. This demonstrates that there exists an execution of \Cref{alg: Proximal Bundle Method} with the modified objective that corresponds exactly to the execution of the algorithm with the original one.
In the dual space, this invariance translates into a shift of the origin, which is theoretically appealing. In \Cref{alg:dual proximal bundle}, the feasible set $\operatorname{conv}(\mathcal{V})$ is replaced by $(a,0) + \operatorname{conv}(\mathcal{V})$ (where $+$ denotes the Minkowski sum), and $\phi$ becomes $\hat{\phi} : w \mapsto \phi(w - a)$.
}

The goal of this study is not to compare the numerical performance of MPB-FA and PBM, but rather to demonstrate that a bundle-type algorithm, which utilizes the same proximal oracles as the classical PBM, can achieve an improved convergence rate in a pervasive setting. It is our hope that the investigation of modified bundle methods, such as MPB-FA, may ultimately lead to an improved complexity analysis of the classical PBM.

\subsection{First convergence rate analysis}
\label{subsection: First convergence rate analysis}
 {From now on, we assume that for Algorithms \ref{alg: Proximal Bundle Method} and \ref{alg:dual proximal bundle}, $x_k=\widehat{x}_k$, $y_{k}=\widehat{y}_k$, and $\cV^k=\widehat\cV^k$ and we use $x_k$, $y_k$, and $\cV^k$ in place of $\widehat{x}_k$,  $\widehat{y}_k$, and $\widehat\cV^k$.}

\subsubsection{Number of serious steps} 
 {We first bound the number of \textit{serious steps}. For convenience, we provide a proof of the following lemma, which appears in \cite{LiangMonteiro2023}.}
\begin{lemma}[Proposition 2 in \cite{LiangMonteiro2023}]
\label{lm: Number of serious steps} 
Under \Cref{assum: general}, \Cref{alg: Proximal Bundle Method} finds an $\epsilon$-optimal solution to  {problem} \eqref{problem: primal} in at most $\|x^*-x_{0}\|^2\rho\epsilon^{-1} +1$ \textit{serious steps},  {where $x^*$ is an optimal solution of  {problem} \eqref{problem: primal}.}
\end{lemma}
\begin{proof}
We consider the iteration $k$ of \Cref{alg: Proximal Bundle Method}, which we suppose is a \textit{serious step}.  {Recall that} $x_{k}$ is the proximal center  {used in \eqref{eq:proxbundle_alg}} to obtain $y_{k+1}$. As iteration $k$ is a \textit{serious step},  {we have} $x_{k+1} = y_{k+1} = -\nabla \cM_{\rho, \phi}(w_{k+1} - \rho x_{k})$.  {Also recall the notation $h:=g+f$ and $h_k:=g+f_k$}.
As $x_{k+1}$ is optimal for \eqref{eq:proxbundle_alg}, whose objective  {$h_k + \frac{\rho}{2}\|\cdot-x_k\|^2$} is $\rho$-strongly convex,
\begin{align}
    h_{k}(x_{k+1}) + \frac{\rho}{2}\|x_{k+1}-x_k\|^2 &\leq h_{k}(x^*) + \frac{\rho}{2}\|x^*-x_k\|^2 -\frac{\rho}{2}\|x_{k+1}-x^*\|^2. \label{ineq: Strong convexity one step}
\end{align}
Rearranging, using $h_k \leq h$ and the \textit{serious step} condition $h(x_{k+1}) - h_k(x_{k+1})\leq \frac{\epsilon}{2}$, and summing over \textit{serious step} index $k_i$: 
\begin{align}
    \min_{i=1,...,T}\{h(x_{k_i}) - h(x^*)\} &\leq \frac{\rho}{2T}\|x_{0}-x^*\|^2 -\frac{\rho}{2T}\sum_{i=0}^{T-1}\|x_{k_{i+1}}-x_{k_i}\|^2 + \frac{\epsilon}{2}. \label{ineq: primal progress sum}
\end{align}
This yields the desired bound on serious steps.
\end{proof}

\subsubsection{Number of null steps}\label{subsec:num null steps} The following lemma bounds the number of consecutive \textit{null steps} between two \textit{serious steps}. During a sequence of  {consecutive} \textit{null steps}, the proximal center is never updated. The minimization step \eqref{eq:proxbundle_alg} in \Cref{alg: Proximal Bundle Method} can be viewed as the update step \eqref{eq:Kelley minimization step} (\Cref{alg: Kelley's algorithm}) applied to $\hat{g}_k:=g_k+\rho/2\|\cdot-x_k\|^2$ where $k$ is the last \textit{serious step}. As $g$ is smooth and  {convex}, $\hat{g}_k$ is smooth and strongly convex. \Cref{theorem: linear convergence of Kelley} shows the linear convergence of Kelley's method applied to $\hat{g}_k + f$. The sequence of \textit{null steps} ends when the \textit{  {serious step test}} $f(y_{k+1}) - f_k(y_{k+1}) \leq \frac{\epsilon}{2}$ is satisfied. We showed in \Cref{lm:bundle_tests} 
that this \textit{  {serious step test}} is exactly the Frank-Wolfe gap of $\tilde{\mathcal{D}}(\cV, x_k)$ in \eqref{eq:proximal bundle dual} (see \Cref{def: FW gap}). We finally upper bound this Frank-Wolfe gap using the primal gap of $\tilde{\mathcal{D}}(\cV, x_k)$. This will show that a \textit{serious step} occurs after at most $\cO(\log\left(\epsilon^{-1}\right))$ \textit{null steps}.

The iteration complexity for \Cref{alg: Kelley's algorithm} given in \Cref{theorem: linear convergence of Kelley} depends on $\|x^*\|$ through the generalized strong convexity constant given in \eqref{eq: generalized strong convexity constant}. The bound on the number of consecutive \textit{null steps} would thus depend on the norm of the optimal solution of the proximal problem $\argmin_{x \in \R^n}\, h(x) +  \frac{\rho}{2}\|x - x_k\|^2$.   {We now provide a uniform bound on this quantity.}
 {\begin{lemma}
\label{lemma: uniform bound on x star k}
     For $i\geq 0$, define $x_i^* := \arg\min_{x \in \R^n}\, h(x) +  \frac{\rho}{2}\|x - x_{j_i}\|^2$, where $x_{j_i}$ are the $i$-th \textit{serious steps} iterates with a global iteration $j_i$ in \Cref{alg: Proximal Bundle Method}. Then,
     \begin{align*}
         \|x_i^*\| \leq \|x^*\| + 2\sqrt{2}\|x^* -x_0\| + 2\sqrt{\frac{\epsilon}{\rho}},
     \end{align*}
     where $x^*$ is an optimal solution of \eqref{problem: primal} which is \(\min_x h(x):=g(x)+f(x)\). 
\end{lemma}}
\begin{proof}
We use \eqref{ineq: Strong convexity one step} for $1\leq k\leq i$ : 
\begin{align*}
    h_{j_k-1}(x_{j_k}) + \frac{\rho}{2} \|x_{j_k} - x_{j_k-1}\|^2 \leq h_{j_k-1}(x^*) + \frac{\rho}{2} (\|x^* - x_{j_k-1}\|^2 - \|x^* - x_{j_k}\|^2)
\end{align*}
This leads to :
\begin{align*}
     &\|x^* - x_{j_k}\|^2 - \|x^* - x_{j_k-1}\|^2 \leq \frac{2}{\rho}(h_{j_k-1}(x^*)-h_{j_k-1}(x_{j_k}) ) - \|x_{j_k} - x_{j_k-1}\|^2\\
    &\quad\leq \frac{2}{\rho}(h_{j_k-1}(x^*) - h_{j_k-1}(x_{j_k}))
    \leq \frac{2}{\rho}( h_{j_k-1}(x^*) + \frac{\epsilon}{2} -h(x_{j_k})) \quad\quad\text{(serious step)}\\
    &\quad\leq \frac{2}{\rho}( h(x^*) + \frac{\epsilon}{2} -h(x_{j_k})) \leq \frac{\epsilon}{\rho}.
\end{align*}
We sum these inequalities from $k=1$ to $i$. We use \Cref{lm: Number of serious steps} to bound the number of \textit{serious steps} by $\frac{\rho \|x_0 - x^*\|^2}{\epsilon} +1$ : 
\begin{align*}
    \|x^* - x_{j_i}\|^2 &\leq \|x^* -x_0\|^2 + \frac{\epsilon}{\rho}\left(\frac{\rho \|x_0 - x^*\|^2}{\epsilon} +1\right)\\
    &= \|x^* -x_0\|^2 + \|x_0 - x^*\|^2+\frac{\epsilon}{\rho} = 2\|x^* -x_0\|^2 + \frac{\epsilon}{\rho}.
\end{align*}
By optimality of $x^*$ and $x_i^*$ for their respective problems: 
\begin{align*}
    h(x^*) + \frac{\rho}{2}\|x^*_i -x_{j_i}\|^2 \leq h(x^*_i) + \frac{\rho}{2}\|x^*_i -x_{j_i}\|^2 \leq h(x^*) + \frac{\rho}{2}\|x^* -x_{j_i}\|^2.
\end{align*}
Thus, $\|x^*_i -x_{j_i}\| \leq \|x^* -x_{j_i}\|$. We combine these two inequalities : 
\begin{align*}
    \|x_i^*\| &\leq \|x^*\| + \|x_i^* - x^*\| \leq  \|x^*\| + \|x_i^* - x_{j_i}\| + \|x^*- x_{j_i}\|
    \leq  \|x^*\| + 2\|x^*- x_{j_i}\|\\
    &\leq  \|x^*\| + 2\sqrt{2\|x^* -x_0\|^2 + \frac{\epsilon}{\rho}}
    \leq  \|x^*\| + 2\sqrt{2}\|x^* -x_0\| + 2\sqrt{\frac{\epsilon}{\rho}}.
\end{align*}
\end{proof}

\begin{lemma}
\label{lemma: number of null steps}
    \Cref{alg: Proximal Bundle Method} performs at most the following number of \textit{null steps} between two \textit{serious steps}
    \begin{align}
            1+\max\left\{2,\frac{D^2 }{ \rho \gamma^2  \bar{\mu}_{\psi, \rho}}\right\} \log \left(\frac{4D^4}{\epsilon^2\rho^2} \right),
    \end{align}
with $\bar{\mu}_{\psi, \rho}$ derived from \Cref{lemma: strong convexity constant of psi} with adapted parameters:
\begin{align*}
    \frac{1}{2}\bar{\mu}_{\psi, \rho}^{-1} &= D_b + 3 \frac{4 M_f^2}{\rho}  + 6M_f (\|x^*\| + 2\sqrt{(1+\rho^2)}\|x^* -x_0\| + 2\sqrt{\rho\epsilon}) \\
    &\quad+ 2(L_g +\rho)\left[\left(\frac{2M_f}{\rho} + \|x^*\| + 2\sqrt{(1+\rho^2)}\|x^* -x_0\| + 2\sqrt{\rho\epsilon}\right)^2 + 1\right].\stepcounter{equation}\tag{\theequation}\label{eq: generalized strong convexity constant for bundle}
\end{align*}
\end{lemma}
\begin{proof}
We consider $k$ the index of a \textit{serious step} and $t\geq k$ an index in the sequence of \textit{null steps} following $k$. Because $\phi$ is convex, we can apply Theorem 6.60 in \cite{BeckMethodsinOptimization} (smoothness of the Moreau envelope) to deduce that $\cM_{\rho, \phi}(\cdot - \rho x_k)$ is $\frac{1}{\rho}$-smooth.  {Recall \(\Psi_k(w,\beta) := \cM_{\rho, \phi}(w - \rho x_k) - \beta\)}. For any $(w_1, \beta_1), (w_2, \beta_2) \in \conv(\cV)$: 
\begin{align*}
    &\|\nabla \Psi_k(w_1, \beta_1)-\nabla \Psi_k(w_2, \beta_2) \| = \left\|\binom{\nabla \cM_{\rho, \phi}(w_1 - \rho x_k) - \nabla \cM_{\rho, \phi}(w_2 - \rho x_k))}{0}\right\|\\
    &\quad\quad= \left\|{\nabla \cM_{\rho, \phi}(w_1 - \rho x_k) - \nabla \cM_{\rho, \phi}(w_2 - \rho x_k))}\right\|
    \leq \frac{1}{\rho}\|w_1 - w_2\|\\
    &\quad\quad\leq \frac{1}{\rho}\sqrt{\|w_1 - w_2\|^2 + (\beta_1 - \beta_2)^2} = \frac{1}{\rho}\|(w_1, \beta_1)^\top  - (w_2, \beta_2)^\top \|.
\end{align*}
We conclude that $\Psi_k$ is $\frac{1}{\rho}$-smooth. Also,  $\cM_{\rho, \phi}(w - \rho x_k) = ({\underbrace{g + \frac{\rho}{2}\|\cdot\|^2}_{(\rho + L_g)-smooth}})^*(\rho x_k  - w).$

The conjugate of a $\sigma$-smooth function is $\frac{1}{\sigma}$-strongly convex. So $\cM_{\rho, \phi}(\cdot - \rho x_k)$ is $(\rho + L_g)^{-1}$-strongly convex. According to \Cref{def: FW gap}, we define the primal gap
\begin{align}
    d^{primal}_{k}(w_{t+1}, \beta_{t+1}) := \Psi_k(w_{t+1}, \beta_{t+1}) - \min_{(w,\beta)\in\conv(\cV^k)}\; \Psi_k(w,\beta), \label{def: primal gap for Psi}
\end{align}
for problem \eqref{problem : null step dual problem} at the point $(w_{t+1},\beta_{t+1})$ of  \eqref{Null step condition is FW Dual gap} and the corresponding Frank-Wolfe gap $d^{FW}_{k}(w_{t+1}, \beta_{t+1})$ given in \eqref{Null step condition is FW Dual gap}. We use Theorem 2 in \cite{lacostejulien2015global} to upper-bound the Frank-Wolfe gap by the  primal gap as 
\begin{align}
    d^{FW}_{k}(w_{t+1}, \beta_{t+1}) &\leq D\sqrt{2\rho^{-1}d^{primal}_{k}(w_{t+1}, \beta_{t+1})} &&\text{(as $\Psi_k$ is $\rho^{-1}$-smooth).} \label{eq: bound on dual gap psi}
\end{align}

 {\Cref{lm:bundle_tests} shows that }the \textit{  {serious step test}} quantity $f(x_{t+1}) - f_t(x_{t+1})$ is exactly this Frank-Wolfe gap. So, if $d^{FW}_{k}(w_{t+1}, \beta_{t+1})\leq \epsilon/2$, then a \textit{serious step} is reached. Using \eqref{eq: bound on dual gap psi}, a $\frac{(\epsilon/2)^2 \rho}{2 D^2}$-optimal solution to \eqref{problem : null step dual problem} would guarantee a Frank-Wolfe gap of at most $\epsilon/2$. 

During a sequence of consecutive \textit{null steps}, the MPB-FA and Kelley's method match,  {because the proximal center is not changed}. \Cref{lm:correspondence Kelley FCFW} shows that these iterates also match with FCFW on problem \eqref{problem : null step dual problem}, with their optimality gap (primal gap) matching as well. 

We define $x_k^* := \underset{x \in \R^n}{\argmin}\, h(x) +  \frac{\rho}{2}\|x - x_k\|^2$. The geometric strong convexity constant given in \Cref{lemma: strong convexity constant of psi} involves $\|x_k^*\|$. To leverage the result of \Cref{theorem: linear convergence of Kelley}, we use the uniform bound on $\|x^*_k\|$ given in \Cref{lemma: uniform bound on x star k}. 

According to \Cref{theorem: linear convergence of Kelley}, Kelley's algorithm finds an $\frac{(\epsilon/2)^2 \rho}{2 D^2}$-optimal solution to \eqref{problem : null step dual problem} in at most the following number of iterations: 
\begin{align*}
    = & 1+\max\left\{2,\frac{D^2 }{ \rho \gamma^2  \bar{\mu}_{\psi, \rho}}\right\} \log \left(\frac{D^4}{(\epsilon/2)^2\rho^2} \right),
\end{align*}
where $\bar{\mu}_{\psi, \rho}$ is the generalized geometric strong convexity constant given by \Cref{lemma: strong convexity constant of psi}, replacing $L_g$ by $\rho + L_g$, $\mu_g$ by $\rho$ and $\|x^*\|$ by $\|x^*\| + 2\sqrt{(1+\rho^2)}\|x^* -x_0\| + 2\sqrt{\rho\epsilon}$. 
\end{proof}

We combine \Cref{lm: Number of serious steps} and \Cref{lemma: number of null steps} to derive \Cref{theorem: bundle complexity 1/epsilon log(epsilon)}. 
\begin{theorem}
    \label{theorem: bundle complexity 1/epsilon log(epsilon)}
    \Cref{alg: Proximal Bundle Method} finds an $\epsilon$-optimal solution to \eqref{problem: primal} in at most the following number of iterations
\begin{align}
    &\left(1+\frac{\|x^*-x_{0}\|^2\rho}{\epsilon}\right) \left[1+\max\left\{2,\frac{D^2}{\bar{\mu}_{\psi, \rho} \rho \gamma^2}\right\} \log \left(\frac{4D^4}{\epsilon^2\rho^2} \right)\right].\label{eq:bundle_first_rate}
\end{align}
\end{theorem}
 {\begin{remark}
The convergence bound $\cO(\epsilon^{-1}\log(\epsilon^{-1}))$ in \eqref{eq:bundle_first_rate} is an order-of-magnitude improvement over the known best convergence rate $\cO(\epsilon^{-2}\log(\epsilon^{-1}))$ of the proximal bundle method \cite{LiangMonteiro2023} applied to solve \eqref{problem: primal}. 
\end{remark}}
\begin{remark}
The analysis of the maximum number of consecutive null steps cannot be directly adapted to a variable parameter $\rho$. However,  \Cref{theorem: bundle complexity 1/epsilon log(epsilon)} still holds if $\rho$ is updated according to a diminishing schedule during \textit{serious steps}.
\end{remark}

\subsection{Improved convergence rate analysis}
\label{subsection: Improved convergence rate}
 {It is surprising that the stronger convergence rate $\cO(\epsilon^{-1}\log(\epsilon^{-1}))$ obtained in \Cref{theorem: bundle complexity 1/epsilon log(epsilon)} for MPB-FA, which is already an order of magnitude better than the best known bound, can be further improved to $\cO(\epsilon^{-8/9}(\log(\epsilon^{-1}))^{2/9})$. This improved bound requires a deeper analysis of the proximal bundle method with novel insights that may be of independent interest. We first state the main theorem, then discuss the key ideas of the proof, and guide the reader through a series of lemmas for the final proof of the improved rate.
\begin{theorem}
    \label{thm: 1/e convergence when g is smooth} 
    Under \Cref{assum: general}, \Cref{alg: Proximal Bundle Method} with \textit{\textcolor{blue}{serious step test}} parameter ${\epsilon}/{2}$ finds an $\epsilon$-optimal solution of \eqref{problem: primal} in at most $\cO({\epsilon^{-8/9}}(\log{\epsilon}^{-1})^{2/9})$ iterations.
\end{theorem}}
\vspace{-5mm}
 {\subsubsection{Key ideas of the improved analysis}\label{subsec:key_improved analysis}
The first analysis \Cref{theorem: bundle complexity 1/epsilon log(epsilon)} bounds the total number of iterations of MPB-FA by the following product
\begin{align}
\text{(\# serious steps)}\times \text{(\# null steps between two serious steps)},\tag{A}\label{blueprint: standard_way}
\end{align}
which implicitly treats all serious steps as non-consecutive. However, a deeper analysis will reveal that not all serious steps are equal. In particular, under certain conditions, serious steps can happen consecutively. When these conditions are not satisfied, we can further distinguish two types of non-consecutive serious steps and upper bound them quite precisely. So, a refined way to bound the iterations of MPB-FA is
\begin{align}
  & \text{(\# \textit{consecutive} serious steps)} + \text{(\# two types of \textit{non-consecutive} serious steps)}\notag \\
    & \hspace{48mm}\times \text{(\# null steps between two serious steps)}.\tag{B}\label{blueprint: improved_way}
\end{align}}
 {The key of this refined analysis is to find quantitative conditions under which two serious steps are consecutive (\Cref{lm: next step is serious step}) and bound the numbers of two types of non-consecutive serious steps when these conditions are violated (\Cref{lm: number of serious steps with distant consecutive iterates} and \Cref{lm: upper bound on the number of serious steps with small gradient}).}
Intuitively, if the consecutive proximal centers $x_{k-1}$ and $x_{k}$ are close, the cuts added to the bundle $\cV^{k}$ while solving $\tilde{\mathcal{D}}(\cV, x_{k-1})$ of \eqref{eq:proximal bundle dual} are relevant for solving $\tilde{\mathcal{D}}(\cV, x_k)$ in the next iteration and this latter problem may be solved in fewer iterations. Determining a threshold on $\|x_{k-1} - x_{k}\|$ to achieve a constant number of iterations for solving $\tilde{\mathcal{D}}(\cV, x_k)$ is not straightforward. A pivotal aspect of the proof lies in considering proximal centers that are close enough, ensuring the optimality of $\tilde{\mathcal{D}}(\cV, x_k)$ after just one step. In that case, $x_{k+1}$ will also be a serious iterate. 
 {\subsubsection{A geometric lemma}\label{subsec:geometric_lemma}
We first show a key lemma, \Cref{lm: minimal angle between FCFW iterates}, which is the foundation of the improved analysis. It studies the geometry of the gradients of the objective function during the iterations of the FCFW (\Cref{alg:FCFW}) applied to solving a convex smooth function over a polytope. This exactly corresponds to applying the modified FCFW (\Cref{alg:dual proximal bundle}) on the Moreau envelope.} 
The main idea of this lemma is that the same  {gradient} direction cannot be explored twice in the linear maximization problem \eqref{FCFW FW step}, as this would indicate that no progress is being made, implying that an optimal solution has already been found. \Cref{lm: minimal angle between FCFW iterates} generalizes this observation by demonstrating that this principle holds even when the two  {gradient} directions are not exactly the same but are sufficiently close to each other. This result leverages the boundedness of the set and the polyhedral nature of the feasible region. 
\begin{lemma}
\label{lm: minimal angle between FCFW iterates}
Let $\psi$ be a convex function and $u_k$ be obtained in step $k$ of FCFW by \eqref{eq:FCFW fully corrective step}, i.e., $u_{k}\in\arg\min_{u\in\conv(\cV^{k-1})}\psi(u)$, and $\tilde{v}_k^{FW}$ be obtained by the corresponding FW step \eqref{FCFW FW step}, i.e. $\tilde{v}_k^{FW}\in\argmax_{v\in\cV}\{d_k^\top v\}$, where $d_{k} = -\nabla \psi(u_{k})/\|\nabla \psi(u_{k})\|$. The following claims hold:
\begin{enumerate}
    \item If $\tilde{v}_{k}^{FW} \in \cV^{k-1}$, then $u_k$ is optimal for the original problem $\min_{u\in  \conv(\cV)}\psi(u)$.
    \item Suppose $d_1 \in \R^n$ is a unit vector and the FW step $\tilde{v}^{FW}_{1} \in {\arg\max}_{v\in \cV}\,\{ d_{1}^\top  v\}$ satisfies $\tilde{v}^{FW}_{1} \in \cV^{k-1}$. Let $\theta$ be the positive angle between $d_{1}$ and $d_k$.  {Let $\gamma$ and $D$ be the pyramidal width and diameter of $\conv(\cV)$, respectively}. Then, if $\theta$ is small enough, i.e. $\sin(\theta)<\frac{\gamma}{3D}$, $u_k$ is again optimal for the original problem.
\end{enumerate}
 
\end{lemma}

\begin{proof}
We first prove claim 1. Suppose that $\tilde{v}_{k}^{FW} \in \cV^{k-1}$. Inequality \eqref{ineq: primal gap < FW gap} and the first-order optimality of $u_k$ in the feasible direction $\tilde{v}_{k}^{FW} - u_k$ of $\cV^{k-1}$ implies 
\begin{align*}
    d^{primal}({u_k}) \leq d^{FW}({u_k}) =\langle d_{k}, \tilde{v}_{k}^{FW} - u_k \rangle \leq  0,
\end{align*}
which shows that $u_k$ is optimal for the original problem $\min_{u\in  \conv(\cV)}\psi(u)$.

\noindent
\begin{minipage}[b]{0.40\textwidth}
To prove claim 2 by the contrapositive, suppose $u_k$ is not optimal, then by claim 1, we know $\tilde{v}_{k}^{FW} \notin \cV^{k-1}$. Suppose $\theta \leq \frac{\pi}{2}$. The definition of the Pyramidal Width and the assumption that $\tilde{v}_{k}^{FW} \notin \cV^{k-1}$ gives $\langle d_{k}, \tilde{v}_{k}^{FW}- u_k \rangle \geq \gamma$. 
\end{minipage}
\hfill
\begin{minipage}[b]{0.58\textwidth}
  \begin{tikzpicture}[>=Stealth, scale=0.75]
    \label{fig:FCFW cannot have small angle}
    \definecolor{darkgreen}{RGB}{0,100,0}
    \definecolor{darkblue}{RGB}{0,0,100}

    \draw[name path=f] (2,3) -- (6.68,1.52);
    \draw[name path=g] (2,3) -- (6.88,3);
    \draw[name path=h] (6.68,1.52) -- (8.16,3);
    \draw[-Stealth, name path=d1] (2.44,1) -- (2.44,2);
    \draw[-Stealth, name path=d2] (6.68,1.52) -- (7.39,2.23);
    \draw[name path=i] (6.68,3.2) -- (6.68,1);
    \draw[name path=j, dashed] (2,3) -- (6.68,1.52);

    \path[name intersections={of=f and g, by=E}];
    \path[name intersections={of=f and h, by=v1}];
    \path[name intersections={of=g and h, by=wk}];
    \path[name intersections={of=i and g, by=s}];

    \coordinate (p) at (7.39,2.23);
    \coordinate[label=below left:$u_1$] (w1) at (2.44,1);
    \coordinate[label=below right:$u_k$] (wk) at (wk);
    \coordinate[label=left:$v_1$] (v1) at (E);
    \coordinate[label=right:$d_1$] (d1) at (2.44,1.5);
    \coordinate[label=right:$d_k$] (dk) at (7.1,1.9);

    \pic[draw,fill=darkblue, angle radius=5mm, angle eccentricity=1.5, "$\theta$"] {angle = p--wk--s};
    \pic[draw,fill=darkgreen, angle radius=10mm, angle eccentricity=1.2, "$\eta$"] {angle = wk--v1--s};
    \pic[draw,fill=white, angle radius=5mm, angle eccentricity=1.4, "$\bar{\eta}$"] {angle = s--wk--v1};

    \foreach \p in {w1,v1,wk} \fill[black] (\p) circle (2pt);
  \end{tikzpicture}
  \raggedleft
  {FCFW : Explored directions}
\end{minipage}%

 \begin{align*}
     \langle d_{1}, \tilde{v}_{1}^{FW}- u_k \rangle \geq \langle d_{1}, \tilde{v}_{k}^{FW}- u_k \rangle
     &= \langle d_{k}, \tilde{v}_{k}^{FW}- u_k \rangle + \langle d_{1} - d_{k}, \tilde{v}_{k}^{FW}- u_k \rangle\\
     &\geq \gamma - \|d_{1} - d_{k}\| \|\tilde{v}_{k}^{FW}- u_k\| 
     \geq \gamma - \|d_{1} - d_{k}\| D. \stepcounter{equation}\tag{\theequation}\label{ineq : d_{1}|v_{1}^{FW}- w_{k}}
 \end{align*}
The first inequality uses the optimality of $\tilde{v}_{1}^{FW}$, the next inequality uses Cauchy-Schwarz and the pyramidal width. The last inequality uses $u_k, \tilde{v}_{k}^{FW} \in \cV$. 

If $u_k = \tilde{v}_{1}^{FW}$, then from inequality \eqref{ineq : d_{1}|v_{1}^{FW}- w_{k}}, $\|d_1 - d_k\| \geq \gamma/D$. We get $\gamma/D \leq \|d_{1} - d_{k}\| = 2 \|\frac{d_{1} + d_{k}}{2} - d_{k}\|= 2\sin\left(\frac{\theta}{2}\right) \leq 2\sin(\theta)$, which proves the result in this case. 

We now suppose that $\|u_k - \tilde{v}_{1}^{FW}\|>0$ and lower bound this quantity. Let $\Pi$ be the plane containing $u_k$ and the directions $d_{1}$ and $d_{k}$. Let $v_{1}$ be the orthogonal projection of $\tilde{v}_{1}^{FW}$ onto $\Pi$. Then, since projections on convex sets are contractive,
\begin{align}
    \|u_k - \tilde{v}_{1}^{FW}\|  \geq \|u_k - v_{1}\|. \label{ineq: projection is contracting}
\end{align}
We define the angle $\eta$ such that $\bar{\eta} := \pi/2 - \eta$ is the positive angle between $d_{1}$ and $v_{1} - u_k$. Since $v_1$ is the orthogonal projection of $\tilde{v}_1^{FW}$ on $\Pi$, we have $\langle d_k, \tilde{v}_1^{FW}-u_k\rangle = \langle d_k, v_1 - u_k\rangle$. Thus, $\langle d_{k}, v_{1} - u_k \rangle = \langle d_{k}, \tilde{v}_{1}^{FW} - u_k \rangle \leq 0$. The inequality follows from the optimality of $u_k\in\argmin_{u\in\conv(\cV^{k-1})}\psi(u)$ and $\tilde{v}_1^{FW}\in\cV^{k-1}$. Remark that $\cos(\bar{\eta} + \theta) \|v_{1} - u_k\| = \langle d_{k}, v_{1} - u_k \rangle \leq 0$. Hence, $ \bar{\eta} + \theta \geq \pi/2$. We conclude that $\eta \leq \theta$.
\begin{align}
    \langle d_{1}, v_{1} - u_k \rangle &= \cos(\frac{\pi}{2} - \eta)\|u_k - v_{1}\| = \sin(\eta)\|u_k - v_{1}\| \leq \sin(\theta)\|u_k - v_{1}\|. \label{ineq: lower bound on sin theta}
\end{align}
Combining \eqref{ineq : d_{1}|v_{1}^{FW}- w_{k}}, \eqref{ineq: projection is contracting},  \eqref{ineq: lower bound on sin theta}, we have
\begin{align*}
     \sin(\theta)\|u_k - \tilde{v}_{1}^{FW}\|&\geq \langle d_{1}, \tilde{v}_{1}^{FW} - u_k \rangle
     \geq (\gamma - \|d_{1} - d_{k}\| D).\stepcounter{equation}\tag{\theequation}\label{ineq : w_{k} - v_{1}^{FW}}
\end{align*}
Now rearrange and apply \eqref{ineq : w_{k} - v_{1}^{FW}} and \eqref{ineq : d_{1}|v_{1}^{FW}- w_{k}}, 
\begin{align}
    \sin(\theta) &\geq \frac{\langle d_{1}, \tilde{v}_{1}^{FW}-u_k  \rangle }{\|u_k - \tilde{v}_{1}^{FW}\|}
    \geq \frac{\gamma}{D} -\frac{D \|d_{1} - d_{k}\|}{\|u_k - \tilde{v}_{1}^{FW}\|} \nonumber\\
    &\geq \frac{\gamma}{D} -\frac{D \|d_{1} - d_{k}\|}{(\gamma - \|d_{1} - d_{k}\| D) \sin(\theta)^{-1}} = \frac{\gamma}{D} -\frac{\sin(\theta) D \|d_{1} - d_{k}\|}{\gamma - \|d_{1} - d_{k}\| D}, \label{ineq: using ineq pyramidal with 2}
\end{align}
where \eqref{ineq: using ineq pyramidal with 2} uses \eqref{ineq : w_{k} - v_{1}^{FW}} a second time. 

Rearranging and using $\|d_{1} - d_{k}\| = 2 \|\frac{d_{1} + d_{k}}{2} - d_{k}\|= 2\sin\left(\frac{\theta}{2}\right) \leq 2\sin(\theta)$,
\begin{align*}
     \sin(\theta) \left(1 + \frac{D \|d_{1} - d_{k}\|}{\gamma - \|d_{1} - d_{k}\| D} \right) &\geq \frac{\gamma}{D},\\
     \implies \; \sin(\theta) \left(1 + \frac{ 2\sin(\theta)D}{\gamma - 2\sin(\theta)D} \right) &\geq \frac{\gamma}{D}.
\end{align*}
If $\gamma - 2\sin(\theta)D \leq  0$, then $\sin(\theta) \geq \frac{\gamma}{D}$. 

Otherwise, 
\begin{align*}
    \gamma \sin(\theta) &\geq \frac{\gamma}{D}(\gamma - 2\sin(\theta)D),\\
    \text{hence,\; }\sin(\theta) &\geq \frac{\gamma}{3D}.
\end{align*}
In both cases, $\sin(\theta) \geq \frac{\gamma}{3D}$. 
\end{proof}



\vspace{-5mm}
 {\subsubsection{A crucial link between proximal centers and Morean envelopes}\label{subsec:proximal_center_Moreaun} Now, we establish in \Cref{lm: upper bound on w_{k} - w_{k+1}} a crucial link between the distance of two proximal centers and the gradients of the associated Moreau envelopes. First, we compute the strong convexity parameter of a function closely related to the Moreau envelope.} 
For $k\in\N$, we define $\hat{s}_k(v,w) := \phi(v) {+} \frac{1}{2 \rho}\|v {+}\rho x_k {-} w\|^2$.  {Recall $\phi=g^*(-\cdot)$.}
Let  $s_k(v,w,\beta) := \hat{s}_k(v,w) {-} \beta$. Define $v_{k+1} := \argmin_{v \in \R^n} s_k(v, w_{k+1}, \beta_{k+1})$. 
Note that $s_k(v_{k+1},w_{k+1},\beta_{k+1}) = \cM_{\rho, \phi} (w_{k+1} - \rho x_{k}) - \beta_{k+1}$.

\begin{lemma}
    \label{lm: strong convexity of phi plus proximal}
     {The function} $\hat{s}_k$ is $\alpha$-strongly convex with
    \begin{align}
        \alpha = \frac{1}{2L_g} + \frac{1}{\rho} -\frac{1}{2} \sqrt{\frac{1}{L_g^2} + \frac{4}{\rho^2}}. 
        \label{eq: alpha strong convexity}
    \end{align}
\end{lemma}
\begin{proof}
As $g$ is $L_g$-smooth, $\phi$ is $L_g^{-1}$-strongly convex. We only need to show that $\frac{1}{2L_g}\|v\|^2 + \frac{1}{2 \rho}\|v +\rho x_k - w\|^2$ is strongly convex  {in $(v,w)$}. We compute $H$, the Hessian of this function, as
\begin{align*}
H &= \frac{1}{L_g}\underbrace{\begin{bmatrix}
 I_n & 0 \\
0 & 0
\end{bmatrix}}_{ {=:P_1}} +  \frac{2}{\rho}\biggl(\underbrace{\frac{1}{2}\begin{bmatrix}
 I_n & - I_n \\
- I_n & I_n
\end{bmatrix}}_{ {=:P_2}}\biggr). 
\end{align*}
A direct algebraic derivation of the smallest eigenvalue of $H$ would prove \eqref{eq: alpha strong convexity}. Instead, we relate this computation to a more general result on the smallest eigenvalue of the non-singular sum of two orthogonal projection matrices. In particular, the Hessian $H$ can be written as $H =aP_1 + bP_2$, where ${a}:=\frac{1}{L_g}$, ${b}:=\frac{2}{\rho}$, and $P_1, P_2$ are defined in the above expression of $H$. Note that $H$ is always non-singular for any positive $L_g$ and $\rho$. The matrix $P_1$ is an orthogonal projection to the subspace $\R^{2n}_n:=\{(u,0)\in\R^{2n} \mid u\in\R^n\}$, while $P_2$ is the orthogonal projection to the subspace $D_{2n}:=\{(u, -u)\in\R^{2n} \mid u \in \mathbb{R}^n\}$.  {Moreover, $\rank(P_1)+\rank(P_2)=2n$ and the first principal angle} between $\R^{2n}_{n}$ and $D_{2n}$ is $\theta_1 = \pi/4$. Thus, the smallest eigenvalue of $aP_1+bP_2$ is \(\frac{1}{2}((a + b) - \sqrt{a^2+b^2})\le\min(a,b)\). This yields \eqref{eq: alpha strong convexity}.
\end{proof}

We consider step $(k-1)$, which we suppose is a \textit{serious step} (so that $x_{k} = y_{k}$). The following lemma establishes the crucial link: when two consecutive serious primal iterates, $x_k$ and $x_{k-1}$, are close, the corresponding gradients associated with the fully corrective step \eqref{dual bundle method fully corrective step} are also close. Our proof draws upon \Cref{lm: strong convexity of phi plus proximal}, the serious nature of step $x_k$, and the inherent smoothness of the Moreau envelope.

\begin{lemma}
\label{lm: upper bound on w_{k} - w_{k+1}}
Assume step $(k-1)$ is a \textit{serious step}. Let $(v_{i+1}, w_{i+1},\beta_{i+1})$ be the optimal solution of \eqref{dual bundle method fully corrective step} in step $i$ with proximal center $x_{i}$, for $i=k-1, k$. Then, 
\begin{align*}
        \| \nabla \cM_{\rho, \phi}(w_{k} - \rho x_{k-1}) - \nabla \cM_{\rho, \phi}(w_{k+1} - \rho x_{k})\|  &\leq (1+\frac{5}{\alpha\rho}) \max(\sqrt{\frac{\epsilon\alpha}{2}},\|x_{k-1} - x_{k}\|).
\end{align*}
\end{lemma}
\begin{proof}
 {From \Cref{lm: dual of prox bundle minimization}, the Moreau envelope \(\cM_{\rho, \phi}(w - \rho x_k)\) is the conjugate of \(\hat{g}_{k}:=g+\rho/2\|\cdot-x_k\|^2\), which is $\rho$-strongly convex. Thus, \(\cM_{\rho, \phi}(w - \rho x_k)\) is $\rho^{-1}$-smooth. Thus, we have
\begin{align*}
\| \nabla \cM_{\rho, \phi}(w_{k} - \rho x_{k-1}) - \nabla \cM_{\rho, \phi}(w_{k+1} - \rho x_{k})\| &\leq \frac{1}{\rho}(\|w_{k} - w_{k+1}\| + \|\rho(x_{k-1} - x_{k})\|).
\end{align*}
To get the inequality in \Cref{lm: upper bound on w_{k} - w_{k+1}}, we only need to upper bound \(\|w_{k} - w_{k+1}\|\) by \(\|x_{k-1} - x_{k}\|\). This is what we manage to do below.} 

By \Cref{lm: strong convexity of phi plus proximal},  {$\hat{s}_k(v,w)$} is $\alpha$-strongly convex  {in $(v,w)$}.  {Since $(v_{k+1}, w_{k+1}, \beta_{k+1})$ is optimal for \eqref{dual bundle method fully corrective step}, 
it is optimal for $\min_{(w,\beta)\in\conv(\cV^k),v}\,\hat{s}_k(v,w)-\beta$, then} we can apply \Cref{lm: Strongly convex in a variable} to $p(v,w) + q(\beta)$ with $p(v,w) = \hat{s}_{k}(v,w)$ and $q(\beta)=-\beta$, leading to the first inequality below
\begin{align}
    &\hspace{5mm}\hat{s}_k(v_{k+1}, w_{k+1}) - \beta_{k+1}\notag\\ 
    &\leq \hat{s}_k(v_{k}, w_{k}) - \beta_{k} - \frac{\alpha}{2}(\|v_{k} - v_{k+1}\|^2 + \|w_{k} - w_{k+1}\|^2) \label{ineq: upper bound on w_{k} - w_{k+1} 0}\\
    & {=\phi(v_k)+\frac{1}{2 \rho}\|v_{k} +\rho x_{k} - w_{k}\|^2-\beta_k- \frac{\alpha}{2}(\|v_{k} - v_{k+1}\|^2 + \|w_{k} - w_{k+1}\|^2)}\notag\\
    &=\phi(v_{k}) + \frac{1}{2 \rho}\|v_{k} +\rho x_{k-1} - w_{k}\|^2 + \langle x_{k} - x_{k-1}, v_{k} +\rho x_{k-1} - w_{k} \rangle \notag\\
    &\hspace{0.5cm} + \frac{\rho}{2}\|x_{k} - x_{k-1}\|^2 - \beta_{k} - \frac{\alpha}{2}(\|v_{k} - v_{k+1}\|^2 + \|w_{k} - w_{k+1}\|^2)\notag\\
    &=s_{k-1}(v_k, w_k, \beta_k) + \langle x_{k} - x_{k-1}, v_{k} +\rho x_{k-1} - w_{k} \rangle + \frac{\rho}{2}\|x_{k} - x_{k-1}\|^2 \label{ineq: upper bound on w_{k} - w_{k+1}}\\
    &\hspace{0.5cm} - \frac{\alpha}{2}(\|v_{k} - v_{k+1}\|^2 + \|w_{k} - w_{k+1}\|^2).\notag
\end{align}
Following the definition of the primal gap \eqref{def: primal gap for Psi}, define
\begin{align*}
d^{primal}_{k-1}(w_k, \beta_k):= s_{k-1}(v_k, w_k, \beta_k) - \min_{(v, w, \beta)\in \R^{n}\times\conv(\cV)}\,s_{k-1}(v,w,\beta).
\end{align*}

As pointed out in  {\Cref{rm:primalgap_FWgap}}, $d^{primal}_{k-1}(w_k, \beta_k) \leq d^{FW}_{k-1}(w_k, \beta_k)$. Step $(k-1)$ is a \textit{serious step} and \Cref{lemma: tests equivalence}  shows that the Frank-Wolfe gap corresponds to the \textit{\textcolor{blue}{serious step test}} in \Cref{alg: Proximal Bundle Method}. We get $d^{primal}_{k-1}(w_k, \beta_k) \leq d^{FW}_{k-1}(w_k, \beta_k) \leq  \frac{\epsilon}{2}$. 
\begin{align*}
    s_{k-1}(v_k, w_k, \beta_k) &= \min_{(v, w, \beta)\in \R^{n}\times\conv(\cV)}\,s_{k-1}(v,w,\beta) + d^{primal}_{k-1}(w_k, \beta_k)\\
    &\leq  \min_{(v, w, \beta)\in \R^{n}\times\conv(\cV)}\,s_{k-1}(v,w,\beta) + \frac{\epsilon}{2}\\
    &\leq s_{k-1}(v_{k+1}, w_{k+1}, \beta_{k+1}) + \frac{\epsilon}{2}. \stepcounter{equation}\tag{\theequation}\label{ineq: relation between primal gap on sk}    
\end{align*}

 {Note \(\hat{s}_k(v_{k+1},w_{k+1})-\beta_{k+1}=s_k(v_{k+1},w_{k+1},\beta_{k+1})\).} Combining \eqref{ineq: upper bound on w_{k} - w_{k+1}} and \eqref{ineq: relation between primal gap on sk}, 
\begin{align*}
    0 &\leq s_{k-1}(v_{k+1}, w_{k+1}, \beta_{k+1}) - s_k(v_{k+1}, w_{k+1}, \beta_{k+1}) + \frac{\epsilon}{2} +  \langle x_{k} {-} x_{k-1}, v_{k} {+}\rho x_{k{-}1} {-} w_{k} \rangle \\ &\quad+ \frac{\rho}{2}\|x_{k} - x_{k-1}\|^2 - \frac{\alpha}{2}(\|v_{k} - v_{k+1}\|^2 + \|w_{k} - w_{k+1}\|^2)\\
    &= \frac{1}{2 \rho}\|v_{k+1} +\rho x_{k-1} - w_{k+1}\|^2- \frac{1}{2 \rho}\|v_{k+1} +\rho x_{k} - w_{k+1}\|^2 + \frac{\epsilon}{2} + \frac{\rho}{2}\|x_{k} - x_{k-1}\|^2 \\
    &\quad+ \langle x_{k} {-} x_{k-1}, v_{k} {+}\rho x_{k{-}1} {-} w_{k} \rangle - \frac{\alpha}{2}(\|v_{k} - v_{k+1}\|^2 + \|w_{k} - w_{k+1}\|^2)\\
    &=\langle x_{k-1} - x_{k}, (v_{k+1} - v_{k})  - (w_{k+1} - w_{k}) \rangle - \frac{\alpha}{2}(\|v_{k} - v_{k+1}\|^2 + \|w_{k} - w_{k+1}\|^2) + \frac{\epsilon}{2}\\
    &\leq \|x_{k-1} - x_{k}\|(\|v_{k+1} - v_{k}\| {+} \|w_{k+1} - w_{k}\|)-\frac{\alpha}{2}(\|v_{k} - v_{k+1}\|^2 {+} \|w_{k} - w_{k+1}\|^2) {+} \frac{\epsilon}{2}\\
     &\leq{-}\frac{\alpha}{4}(\|v_{k} - v_{k+1}\| {+} \|w_{k} - w_{k+1}\|)^2 + \|x_{k-1} - x_{k}\|(\|v_{k+1} {-} v_{k}\| {+} \|w_{k+1} {-} w_{k}\|){+} \frac{\epsilon}{2}.\stepcounter{equation}\tag{\theequation}\label{ineq: degree 2 polynomial ineq}
\end{align*}
The penultimate inequality follows from the Cauchy-Schwarz inequality, while the final inequality leverages the identity $2(a^2 + b^2) \geq (a + b)^2$.

The expression \eqref{ineq: degree 2 polynomial ineq}  {defines} a quadratic polynomial with a negative leading coefficient  {\[P(z)=-\frac{\alpha}{4} z^2 + \|x_{k-1}-x_k\|z+\frac{\epsilon}{2},\] which satisfies} 
$P(0) {=\epsilon/2} > 0$.  {Since $\hat{z}:=\|v_{k} - v_{k+1}\| + \|w_{k} - w_{k+1}\| \geq 0$, $P(\hat{z})\ge 0$ implies} $\hat{z}\leq x^+$, where $x^+$ is the largest root of $P$  {on the right-hand side of \eqref{eq:x+},} 
\begin{align}
    \|v_{k} - v_{k+1}\| + \|w_{k} - w_{k+1}\| &\leq 2\frac{\|x_{k-1} - x_{k}\| + \sqrt{\|x_{k-1} - x_{k}\|^2 +  \epsilon \alpha/2}}{\alpha}\label{eq:x+}\\
    &\leq \frac{2(1+\sqrt{2})}{\alpha}\max(\sqrt{\frac{\epsilon\alpha}{2}},\|x_{k-1} - x_{k}\|)\notag.
\end{align}
As $ \|v_{k} - v_{k+1}\| \geq 0$, $\|w_{k} - w_{k+1}\| \leq \frac{5}{\alpha}\max(\sqrt{\frac{\epsilon\alpha}{2}},\|x_{k-1} - x_{k}\|).$

 {Now, we can go back to use} the $\frac{1}{\rho}$-smoothness of $\cM_{\rho, \phi}$:
\begin{align*}
    \| \nabla \cM_{\rho, \phi}(w_{k} - \rho x_{k-1}) - \nabla \cM_{\rho, \phi}(w_{k+1} - \rho x_{k})\| &\leq \frac{1}{\rho}(\|w_{k} - w_{k+1}\| + \|\rho(x_{k-1} - x_{k})\|)\\
    &\leq (1+\frac{5}{\alpha\rho}) \max(\sqrt{\frac{\epsilon\alpha}{2}},\|x_{k-1} - x_{k}\|).
\end{align*}
\end{proof}

 {\subsubsection{A sufficient condition for consecutive serious steps}\label{subsec:condition_consecutive_serious steps} Now we can establish a sufficient condition under which serious steps become consecutive.}
The following lemma gives conditions that ensure that the angle between the consecutive iterates $y_k$ and $y_{k+1}$ is small enough to apply \Cref{lm: minimal angle between FCFW iterates} and conclude that step $k$ is a \textit{serious step}. 
\begin{lemma}[Sufficient condition for consecutive serious steps]
\label{lm: next step is serious step}
    Let $A \geq \sqrt{\alpha/2}$  {with $\alpha$ defined in \eqref{eq: alpha strong convexity}}. In the execution of \Cref{alg: Proximal Bundle Method}, suppose, step $(k-1)$ is a \textit{serious step}, $\|y_{k}\|\geq \sqrt{\epsilon}   \frac{6A(1+5(\alpha\rho)^{-1}) D}{\gamma}$, and $\|x_{k-1} - x_{k}\| \leq A\sqrt{\epsilon}$,  {where $\gamma$ and $D$ are the pyramidal width and diameter of $\conv(\cV)$}. Then, the step $k$ is also a \textit{serious step}.  
\end{lemma}
\begin{proof}
Following \eqref{dual bundle method  first order oracle} in \Cref{alg:dual proximal bundle},
\begin{align*}
    y_{k} &= -\nabla \cM_{\rho, \phi}(w_{k} - \rho x_{k-1}) \quad\text{and}\quad
    y_{k+1} = -\nabla \cM_{\rho, \phi}(w_{k+1} - \rho x_{k}).
\end{align*}

Let $\theta$ be the angle between $y_{k}$ and $y_{k+1}$,
\begin{align}\|y_{k} - y_{k+1}\|^2 = \|y_{k}\|^2 \sin(\theta)^2 + (\|y_{k+1}\| - \|y_k\|\cos(\theta))^2 \geq \|y_{k}\|^2 \sin(\theta)^2. \label{ineq: sin theta vs y}\end{align}

Using \eqref{ineq: sin theta vs y},  {the assumptions in the lemma statement}, and \Cref{lm: upper bound on w_{k} - w_{k+1}},
\begin{align}
    \sin(\theta) &\leq \frac{\|y_{k} - y_{k+1}\|}{\|y_{k}\|}
    \leq \frac{\gamma (1+5(\alpha\rho)^{-1}) \max(\|x_{j} - x_{j+1}\|,\sqrt{\epsilon \alpha/2} )}{6(1+5(\alpha\rho)^{-1}) \sqrt{\epsilon} A  D}\nonumber\\
    &= \frac{\gamma \max(\|x_{k-1} - x_{k}\|,\sqrt{ \epsilon\alpha/2} )}{6 \sqrt{\epsilon} A  D}
    \leq \frac{\gamma \max(1,\frac{\sqrt{ \alpha/2}}{A})}{6   D}
    \leq \frac{\gamma}{6D}. \label{ineq: inequality on angle in n dim}
\end{align}

As ${\gamma}/({6D})<1$ by the definition of pyramidal width, Eq. \eqref{ineq: inequality on angle in n dim} leads to $\|y_{k} - y_{k+1}\|< \|y_{k}\|$, which shows that $\cos(\theta) > 0$. Thus $0\leq \theta <\pi/2$.  The inequality \eqref{ineq: inequality on angle in n dim} upper bounds the angle between $y_k$ and $y_{k+1}$. At step $(k-1)$,  {the linear optimization oracle \eqref{eq: dual bundle method FW step} uses} the negative gradient  {of the objective of \eqref{dual bundle method fully corrective step}}: $\binom{-\nabla \cM_{\rho, \phi}(w_{k} - \rho x_{k-1})}{-1} = \binom{y_k}{-1}$. To fulfill the hypothesis of \Cref{lm: minimal angle between FCFW iterates}, we will upper bound the angle $\theta'$ between $\binom{y_{k}}{-1}$ and $\binom{y_{k+1}}{-1}$ by taking into account the $(n+1)^{th}$ dimension of this gradient. We can restrict ourselves to a 3-dimensional subspace of $\R^{n+1}$ spanned by $\binom{y_{k}}{0}$, $\binom{y_{k+1}}{0}$, and $\binom{{0}_{n}}{1}$. We denote $\tilde{y}_k$ and $\tilde{y}_{k+1}$ the restriction of $\binom{y_{k}}{-1}$ and $\binom{ y_{k+1}}{-1}$ in this subspace.  {Define $r:=\|y_k\|
$ and $r':=\|y_{k+1}\|$}. We can find an orthonormal basis of this subspace such that the coordinates of $\tilde{y}_k$ and $\tilde{y}_{k+1}$ are:
\begin{align*}
   \tilde{y}_k = \begin{bmatrix}
    r \\
    0 \\
    -1 \\
\end{bmatrix}
\quad \text{and} \quad
\tilde{y}_{k+1} = \begin{bmatrix}
    r'\cos(\theta) \\
    r'\sin(\theta) \\
    -1 \\
\end{bmatrix}. \quad \text{Note that} \quad \tilde{y}_k \times \tilde{y}_{k+1} =
\begin{bmatrix}
r'\sin(\theta) \\
r - r'\cos(\theta) \\
rr'\sin(\theta) \\
\end{bmatrix},
\end{align*}
where $\tilde{y}_k \times \tilde{y}_{k+1}$ denotes the cross product of $\tilde{y}_k$ and $\tilde{y}_{k+1}$.
\begin{align*}
\|\tilde{y}_k \times \tilde{y}_{k+1}\|^2 &= r'^2 + r^2 - 2rr'\cos(\theta) + r^2r'^2\sin^2(\theta) \\
\implies\sin^2(\theta')  {=\frac{\|\tilde{y}_k \times \tilde{y}_{k+1}\|^2}{\|\tilde{y}_k\|^2 \|\tilde{y}_{k+1}\|^2}} &=  \frac{r'^2 + r^2 - 2rr'\cos(\theta) + r^2r'^2\sin^2(\theta)}{\|\tilde{y}_k\|^2 \|\tilde{y}_{k+1}\|^2}\\
 &= \frac{(r - r')^2 + 4rr'\sin^2\left(\frac{\theta}{2}\right) + r^2r'^2\sin^2(\theta)}{(r^2 + 1)(r'^2 + 1)}\tag{\theequation}\label{ineq: cos half angle}\\
&< \frac{(r - r')^2 }{(r^2 + 1)\cdot 1}+ \frac{4rr'}{r^2 + r'^2}\sin^2\left(\theta\right) + \sin^2(\theta)\\
&\leq \frac{\|y_{k} - y_{k+1}\|^2 }{\|y_{k}\|^2}+ 2\sin^2\left(\theta\right) + \sin^2(\theta) \stepcounter{equation}\tag{\theequation}\label{ineq: r minus r prime}\\
&\leq \left(\frac{\gamma}{6D}\right)^2 + 3\sin^2(\theta) \leq 4\left(\frac{\gamma}{6D}\right)^2.
\end{align*}
Equation \eqref{ineq: cos half angle} uses $\cos(\theta) = 1-2\sin(\theta/2)^2$ and \eqref{ineq: r minus r prime} uses $|r - r'| \leq \|y_{k} - y_{k+1}\|$. Taking the square root, we conclude that $\sin(\theta')< \frac{\gamma}{3 D}$ with $0\leq\theta < \frac{\pi}{2}$.

 We consider the execution of the FCFW algorithm applied to the problem \eqref{problem : null step dual problem}. At step $k$, when the fully corrective step \eqref{dual bundle method fully corrective step} is computed, the direction $\tilde{y}_k$ has already been explored at the previous iteration $k$. We have shown that $\theta'$ the angle between $\tilde{y}_{k+1}$, the gradient of the iterate $(w_{k+1}, \beta_{k+1})$, and $\tilde{y}_k$ is such that $\sin(\theta')<\frac{\gamma}{3 D}$ and $\theta'\leq \frac{\pi}{2}$. We apply \Cref{lm: minimal angle between FCFW iterates} with $d_1 = \tilde{y}_k$ and the iterate $u_k = (w_{k+1}, \beta_{k+1})$ to deduce that $u_k$ is optimal for the \textit{null steps} problem \eqref{problem : null step dual problem}. In particular, this shows the Frank-Wolfe gap is zero and that $k$ is a \textit{serious step}. 
\end{proof}
 {\subsubsection{Bounding two types of non-consecutive serious steps}\label{subsec:nonconsecutive_serious steps}  
We now consider non-consecutive \textit{serious steps}, i.e., serious steps that are} followed by at least one \textit{null step}. We can bound the number of non-consecutive serious by bounding the number of \textit{serious steps} for which one of the hypotheses of  \Cref{lm: next step is serious step} is violated.  {There are two types.} Intuitively, looking at \eqref{ineq: primal progress sum}, a large distance between consecutive \textit{serious steps} iterates ensures a sufficient decrease in function value. Also, when a small gradient $y_k$ is encountered, the number of remaining \textit{serious steps} is small. It is the subject of the next two lemmas. We use step $k_i$ to denote the $i^{th}$ \textit{serious step}.  
\begin{lemma}[ {Type I: Distant proximal centers}]
\label{lm: number of serious steps with distant consecutive iterates}
    There are at most $\frac{2 \|x_0 - x^*\|^2}{\epsilon A^2} + \frac{2}{A^2\rho}$ serious iterates such that $\|x_{k-1} - x_{k}\| > A \sqrt{\epsilon}$ before reaching an $\epsilon$-optimal solution. 
\end{lemma}

\begin{proof}
Using \eqref{ineq: primal progress sum}, for all $T\in \N\setminus{\{0\}}$ : 
\begin{align}
    0\leq \min_{i=1,...,T}\{h(x_{k_i}) - h(x^*)\} &\leq \frac{\rho}{2T}\|x_{0}-x^*\|^2 -\frac{\rho}{2T}\sum_{i=1}^{T}\|x_{k_{i}}-x_{k_{i}-1}\|^2 + \frac{\epsilon}{2}. \nonumber\\
   \implies \sum_{i=1}^{T}\|x_{k_{i}}-x_{k_{i}-1}\|^2 &\leq \frac{T\epsilon}{\rho} + \|x^*-x_{0}\|^2.\label{ineq: number of large serious steps}
\end{align}
According to \Cref{lm: Number of serious steps}, we reach an $\epsilon$-optimal solution in at most $T =  \frac{\rho \|x_0 - x^*\|^2}{\epsilon} +1 $ \textit{serious steps}. Let $K\in \N$ be the number of iterates such that $\|x_{k-1} - x_{k}\| > A \sqrt{\epsilon}$. Following \eqref{ineq: number of large serious steps}, $K A^2 \epsilon \leq  \frac{T\epsilon}{\rho} + \|x^*-x_{0}\|^2$. We deduce that 
\begin{align*}
    K&\leq \frac{T}{A^2\rho} + \frac{\|x^*-x_{0}\|^2}{A^2\epsilon}
    \leq \frac{2}{A^2\rho}\left(\frac{\rho \|x_0 - x^*\|^2}{\epsilon} + 1\right)
    = \frac{2 \|x_0 - x^*\|^2}{\epsilon A^2} + \frac{2}{A^2\rho}.
\end{align*}
\end{proof}

\begin{lemma}[ {Type II: Small gradient of Moreau envelope}]
\label{lm: upper bound on the number of serious steps with small gradient}
Let $M>0$. In the execution of \Cref{alg: Proximal Bundle Method}, suppose that the serious iterate $y_{k}$ is such that $\|y_{k}\|\leq \sqrt{\epsilon} M$, we find an $\epsilon$-optimal solution after at most $\lceil\rho M^2\rceil$ \textit{serious steps}. 
\end{lemma}

\begin{proof}
We recall the convergence guarantee in terms of the number of \textit{serious steps} given in \eqref{ineq: primal progress sum}. Suppose the iterate $y_k = x_{k_i}$ corresponds to the $i^{th}$ \textit{serious step}. For \textit{serious steps} $x_{k_l}$ and $x_{k_{l-1}}$, it holds
\begin{align*}
    h(x_{k_l})- h(x^*)&\leq  \frac{\rho}{2}(\|x^*-x_{k_{l-1}}\|^2 -\|x_{k_l}-x^*\|^2) +  \frac{\epsilon}{2} - \frac{\rho}{2}\|x_{k_{l-1}} - x_{k_l}\|^2.
\end{align*}
Let $T_i\geq 1$. We have
\begin{align*}
    &\frac{1}{T_i}\sum_{l=i}^{T_i+i-1} h(x_{k_l})- h(x^*)\leq  \frac{\rho}{2T_i}\|x^*-x_{k_i}\|^2 +  \frac{\epsilon}{2}  - \frac{1}{T_i}\sum_{l=i}^{T_i+i-1} \frac{\rho}{2}\|x_{k_{l-1}} - x_{k_l}\|^2,\\
   \implies &\underset{l \in [1,T_i+i-1]}{\min}\,\{h(x_{k_l})\}- h(x^*) \leq \frac{\|x^*-x_{k_i}\|^2\rho}{2T_i} + \frac{\epsilon}{2}.
\end{align*}
An $\epsilon$-optimal solution is found if $\frac{\|x^*-x_{k_i}\|^2\rho}{2T_i} \leq \frac{\epsilon}{2}$. The algorithm stops after at most $T_i$ iterations with 
\begin{align*}
    \lfloor T_{i} \rfloor&\leq \frac{\|x^*-x_{k_i}\|^2\rho}{\epsilon} = 
    \frac{\|y_{k}\|^2\rho}{\epsilon}
    \leq \rho \left(\sqrt{\epsilon} M\right)^2 \epsilon^{-1}
    =  \rho M^2.
\end{align*}
This inequality holds as we supposed that $x^* = 0$. We can upper bound the number of remaining iterations by $\lceil \rho M^2 \rceil$, concluding the proof of the lemma.
\end{proof}

 {\subsubsection{Proof of the main theorem}
We now have all the ingredients to prove} the main result, \Cref{thm: 1/e convergence when g is smooth}, which improves upon the convergence rate given in \Cref{theorem: bundle complexity 1/epsilon log(epsilon)}.  {As outlined in \Cref{subsec:key_improved analysis}}, we distinguish between ``consecutive'' \textit{serious steps} that satisfy the assumptions of \Cref{lm: next step is serious step}, and  {non-consecutive serious steps of Type I} ``distant centers''  {(\Cref{lm: number of serious steps with distant consecutive iterates})} or  {Type II} ``small gradient''  {(\Cref{lm: upper bound on the number of serious steps with small gradient})} that will be followed by $\cO\left(\log({\epsilon}^{-1}) \right)$ \textit{null steps} (\Cref{lemma: number of null steps}). We find an optimal value for the constant A (which does not appear in the algorithm). We also give the optimal parameter $\rho$ and the corresponding convergence rate.

\input{proof_final_improved_rate_v3}


 

\section{Numerical Experiments for MPB-FA} \label{sec: bundle management}
In this section, we present some preliminary  {computational} results. The goal is not to numerically compare MPB-FA to other proximal bundle methods, but to compare three bundle management policies for MPB-FA. The three policies mainly differ in how to handle the bundle of cuts after \textit{serious steps}. We provide preliminary numerical results to illustrate that the \textit{active-cut policy} that maintains a proper level of model accuracy by keeping \textit{active cuts} after \textit{both} \textit{serious} and \textit{null steps} has advantages over the \textit{all-cut policy} that may keep too many cuts or the \textit{single-cut policy} that may keep too few cuts. This result aligns with the improved convergence rate from \Cref{thm: 1/e convergence when g is smooth} which requires keeping active cuts in the bundle after a \textit{serious step} as well. 

In \Cref{alg: Kelley's algorithm} and \Cref{alg: Proximal Bundle Method}, a cut $i \leq k$ is said to be \textit{active} at iteration $k$ if $f_k(y_{k+1}) = (v_i^{FW})^\top  y_{k+1} + b_i^{FW}$. These active cuts correspond to the supporting hyperplanes of $f$ that have the maximal value at $y_{k+1}$ among all cuts defining $f_k$.

\vspace{-5mm}
 {\subsection{Bundle management policies}\label{subsec:bundle_management}
In this subsection, we compare three policies for bundle management.}

\underline{All-cut policy}: retaining all the past cuts, i.e. both active and non-active cuts of the model $f_k$ after every \textit{serious} and \textit{null step} as in the current form of \Cref{alg: Proximal Bundle Method}. This policy may lead to growing memory demands and escalate the difficulty in solving the associated minimization problems \eqref{eq:proxbundle_alg}. 

\underline{Single-cut policy}: 
as noted in \cite{Diaz2023optimal, LiangMonteiro2023}, retaining some approximation accuracy of the previous model is only necessary after \textit{null steps}. To keep the bundle as lean as possible, it is suggested in \cite{liang2021proximal} to drop all but the \textit{newest added cut} after every \textit{{serious step}} (thus the name) and keep active cuts only after \textit{{null steps}}. 

\underline{Active-cut policy}: retaining active cuts after \textit{both} \textit{serious} and \textit{null steps}.  Since at most $n+1$ cuts are retained during \textit{null step} sequences \cite{Du2017,zhou2018limited, VANACKOOIJ2017}, retaining active cuts after a \textit{serious step} does not increase the maximum memory requirement.
\vspace{-1mm}
\textcolor{blue}{
\begin{lemma}\label{lm:active-cut}
    The convergence rate of \Cref{thm: 1/e convergence when g is smooth} holds for the active-cut policy.
\end{lemma}
\begin{proof}
The \emph{active-cut policy} retains only the subset of cuts that remain active, 
corresponding to the face of $\mathcal{V}$ on which $w_k$ lies. This restriction is sufficient to ensure that the convergence rates established 
for the FCFW algorithm in~\cite{braun2023conditional} remain valid. It therefore remains to verify that \Cref{lm: minimal angle between FCFW iterates} still applies. Indeed, the requirement in~\eqref{ineq: upper bound on w_{k} - w_{k+1} 0} is satisfied 
because $w_k$ lies within the convex hull of $\mathcal{V}^k$.
\end{proof}}
 To the best of our knowledge, our analysis is the first to give a theoretical ground for retaining active cuts after a \textit{serious step}. Indeed, the all-cut and single-cut policies prove to be less balanced and less effective than the active-cut policy as illustrated by the following numerical experiment.
\begin{figure}[h!]
\centering
 \begin{subfigure}{0.325\textwidth}
  \centering\includegraphics[width=\linewidth]{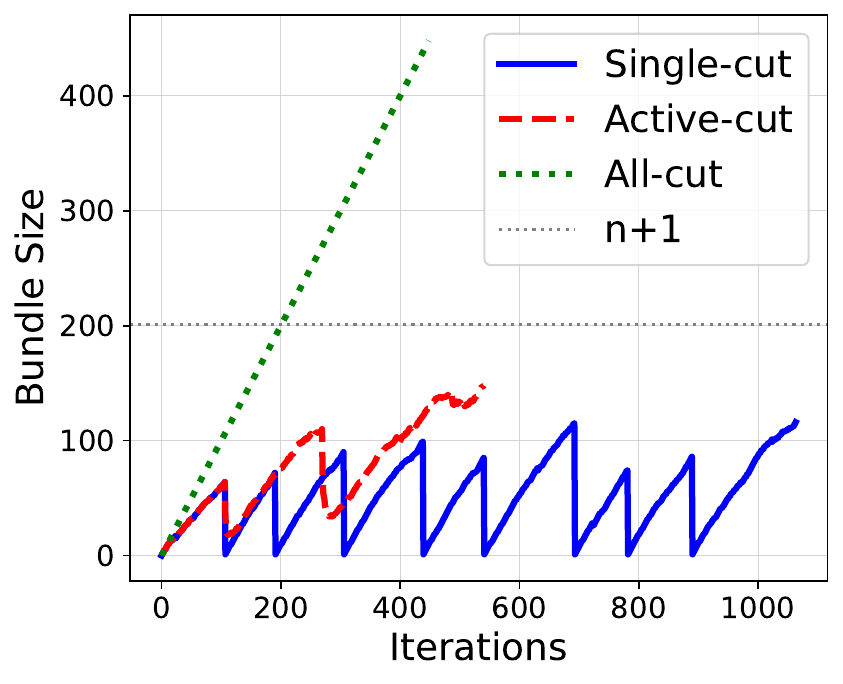}
 \end{subfigure}
 \begin{subfigure}{0.325\textwidth}
   \centering\includegraphics[width=\linewidth]{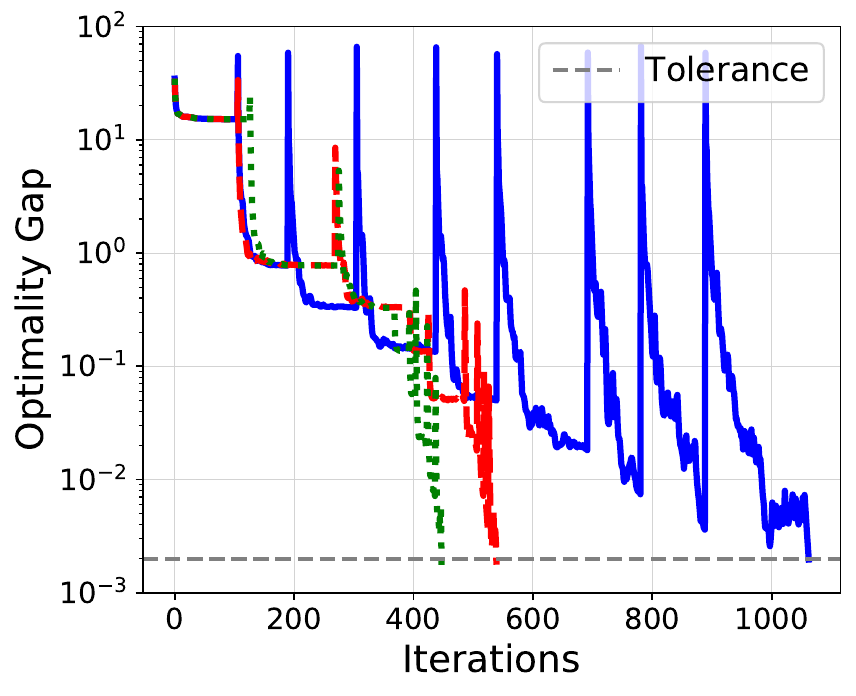}
 \end{subfigure}
 \begin{subfigure}{0.325\textwidth}
   \centering\includegraphics[width=\linewidth]{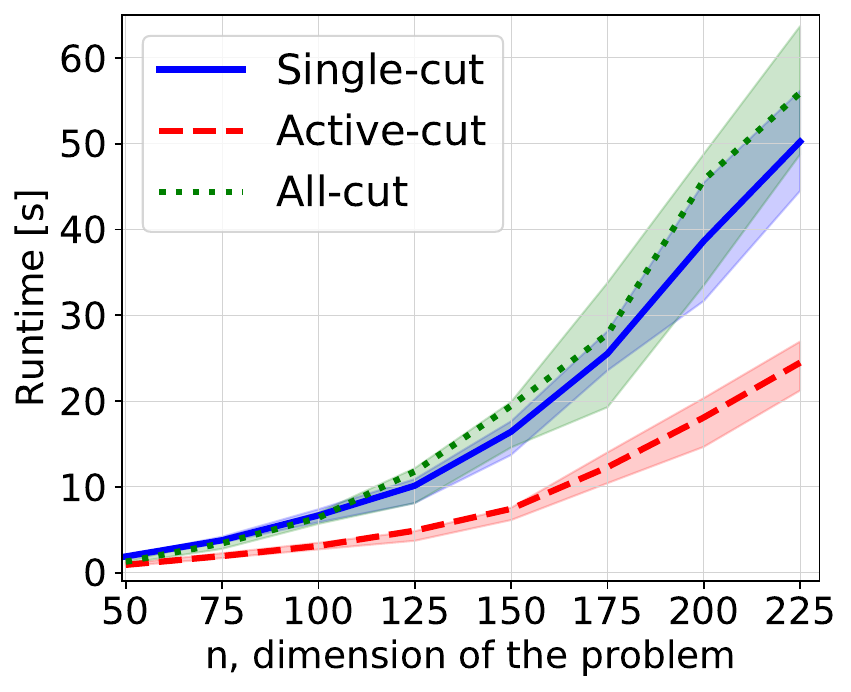}
 \end{subfigure}
 \vspace{-2mm}
 \caption{\small The two leftmost plots compare the memory usage and the optimality gap across bundle management policies ($n = 200$). The last plot compares their runtime for varying problem dimension $n$ with interquartiles.}
  \label{fig: experiment}
\end{figure}

\emph{Experiment setup} 
We illustrate the performance gain of the active-cut policy through a numerical study on the problem: 
\begin{align*}
    \min_x\, \biggl\{\frac{1}{2}\|x\|^2 + \max_y\{x^\top  y + d : y\in [-1, 1]^n, d\in [-1, 1],  A y + c d\leq b\}\biggr\}.
\end{align*} 
The dimension of the problem, $n$, is five times the number of constraints $m$, and $A$, $b$, and $c$ are chosen with uniformly random components in $[-1, 1]$. The problem is solved with accuracy $\epsilon = 2\cdot 10^{-3}$, while $\rho$ is set to $1$. The implementation \footnote{\href{https://github.com/dfersztand/Improved-Complexity-Analysis-for-the-Proximal-Bundle-Algorithm-Under-a-Novel-Perspective}{https://github.com/dfersztand/Improved-Complexity-Analysis-for-the-Proximal-Bundle-Algorithm-Under-a-Novel-Perspective}} is done in Julia. All experiments were conducted on an Intel Xeon Platinum processor machine featuring 48 cores and 4 GB of RAM per core. This problem captures the structure of many composite nonsmooth problems, such as SVM and robust optimization. 

\Cref{fig: experiment} illustrates this tradeoff by comparing the memory usage (bundle size), optimality gap, and total runtime of the three policies. For the two leftmost plots of \Cref{fig: experiment}, which illustrate the memory requirement and iteration complexity across the three bundle management policies, we have chosen $n = 200$ and $m = 40$. The rightmost plot shows the runtime of the MPB-FA for the three bundle management policies with instances of increasing problem size. The runtimes have been averaged over 20 random instances, with the shaded area representing the interquartile range.

Retaining more cuts reduces the total number of iterations but increases the time per iteration. The leftmost figure shows the fast growth of bundle size of the all-cut policy (green), the frequent cut cleansing in the single-cut policy (blue), and balanced bundle size of active-cut policy (red). In the middle figure, the all-cut policy has the fewest number of iterations, while cleansing the bundle after \textit{serious steps} as in the single-cut policy leads to wild increase of optimality gap after \textit{serious steps} and slow convergence. The right figure shows that fewer iterations in the all-cut policy does not result in a shorter runtime because the iterations become more computationally expensive. Retaining only the active cuts at \textit{serious step} achieves a good balance and results in the best runtime.
\vspace{-2mm}
 {\subsection{Empirical behavior of MPB-FA and classical PBM supporting the improved complexity analysis}\label{subsec:prob of consec serious steps}
Preliminary experiments show that the theoretical insights of \Cref{lm: next step is serious step} are indeed supported by the empirical performance of the MPB-FA (\Cref{alg: Proximal Bundle Method}) under the same experimental setup as in \Cref{subsec:bundle_management}.}

\begin{figure}[h!]
\centering
 \begin{subfigure}{0.485\textwidth}
  \centering\includegraphics[width=\linewidth]{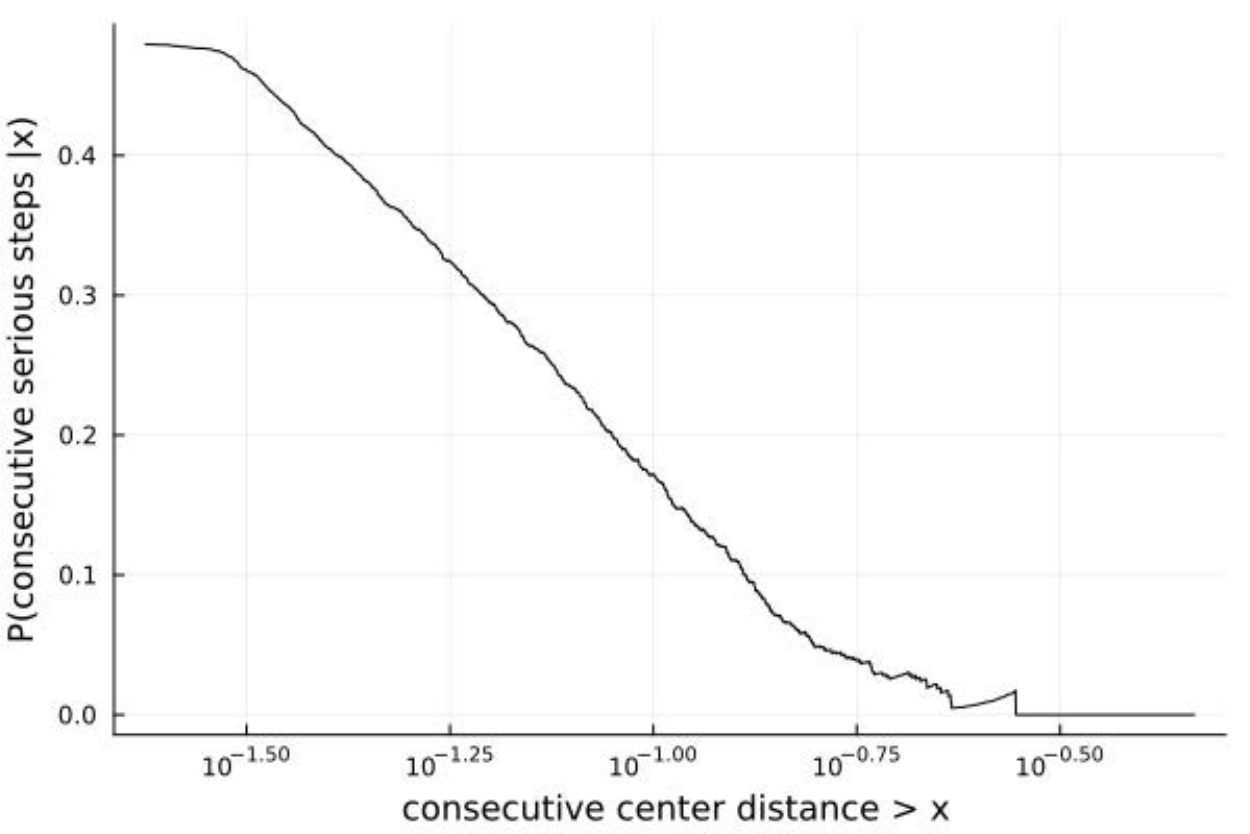}
  \subcaption{MPB-FA}\label{fig:MPB-FA}
 \end{subfigure}
 \begin{subfigure}{0.485\textwidth}
   \centering\includegraphics[width=\linewidth]{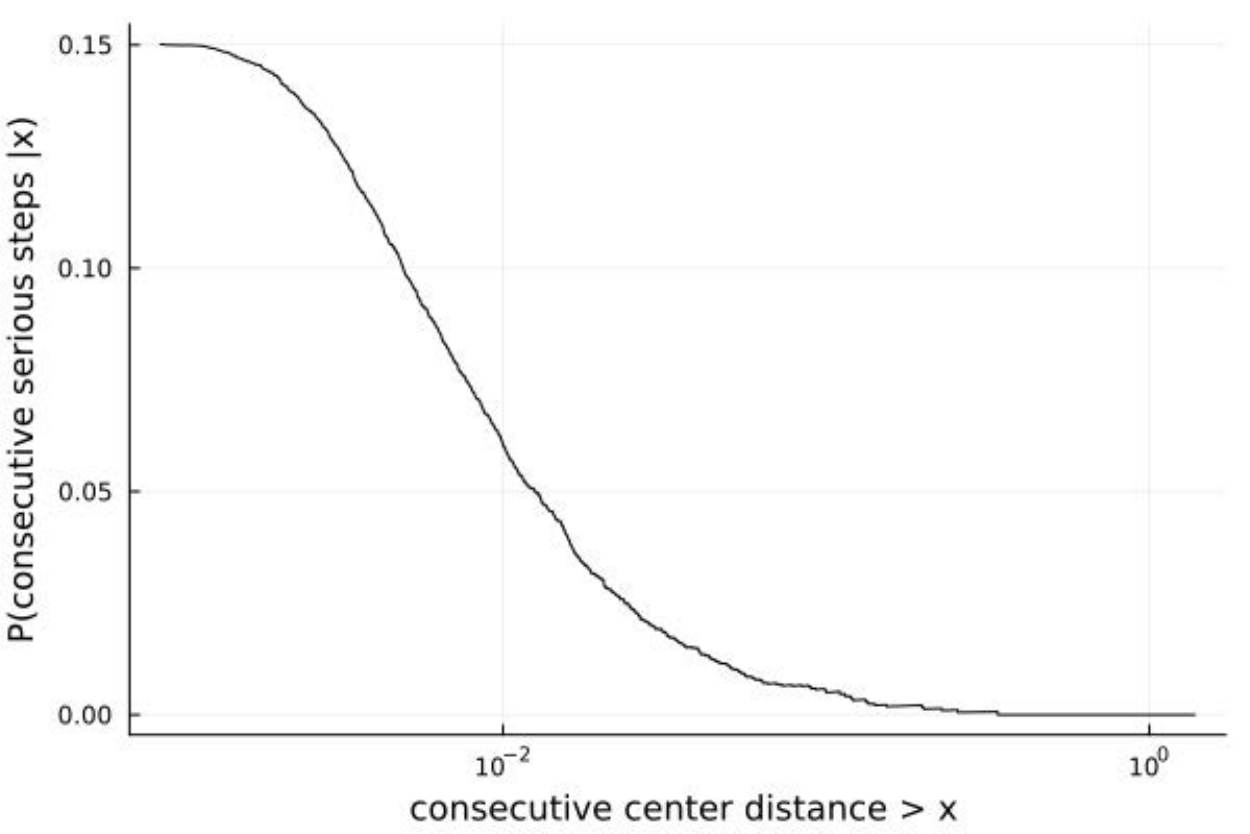}
   \subcaption{Classical PBM}\label{fig:classicalPBM}
 \end{subfigure}
 \vspace{-2mm}
 \caption{\small \textcolor{blue}{Empirical probability of observing consecutive serious steps with respect to the distance between the two preceding proximal centers.}}
  \label{fig: prob of consecutive centers experiment}
\end{figure}
 {\Cref{fig:MPB-FA} plots the empirical probability of consecutive serious steps as a function of the distance between the two previous proximal centers. 
\Cref{fig:MPB-FA} shows a clear trend of a higher frequency of consecutive serious steps in MPB-FA when the two previous prox centers are close. We also conduct the same experiments on classical PBM, which differs from MPB-FA by using the relative-error serious step test \eqref{eq:classicPBA_test}. For classical PBM, \Cref{fig:classicalPBM} shows a similar behavior of a higher frequency of consecutive serious steps when the two preceding proximal centers are close. This suggests that the insights for the improved analysis in \Cref{subsection: Improved convergence rate} may have revealed some important behavior of the classical PBM that is hidden from the existing analysis on PBM.}

\section{Discussion and conclusions}
\subsection{Complexity lower bounds} \Cref{alg:FCFW} and \Cref{alg:dual proximal bundle} (and their dual algorithms \Cref{alg: Kelley's algorithm} and \Cref{alg: Proximal Bundle Method}) are LP-based convex programming (LCP) methods as defined in \cite{lan2014complexity}. \Cref{alg:FCFW} is an LCP that minimizes a $1/\mu_g$-smooth function. It is shown in \cite{lan2014complexity} that the worst-case iteration complexity of LCP to solve smooth (even strongly-convex) problems cannot be smaller than $\min\left\{{n}/{2}, {D^2}\epsilon^{-1}/(4\mu_g)\right\}$. This aligns with the linear convergence rate given in \Cref{theorem: linear convergence of Kelley}. Indeed, the constants in this rate involve $\frac{D^2}{\bar{\mu}_\psi\mu_g \gamma^2} \geq \frac{D^2}{\gamma^2} \geq \lfloor\frac{n}{2}\rfloor$, with the final inequality regarding the condition number \({D^2}/{\gamma^2}\) being established in \cite{Alexander1977TheWA}. Likewise, the convergence rate \eqref{final convergence rate last thm} presented in \Cref{thm: 1/e convergence when g is smooth} includes a term $163 D^2/\gamma^2\geq n/2$.

\subsection{Best known convergence rates}\label{subsec:best known rates}
In the setting of \Cref{section : Proximal Bundle algorithm method rate}, the analysis of Liang and Monteiro \cite{LiangMonteiro2023} provides an $\cO(\epsilon^{-2}(\log\epsilon^{-1}))$ iteration complexity. It is worth noting that in the regime where 
$n \gg {\epsilon}^{-1}$ (which corresponds to high-dimensional settings or low-accuracy scenarios), the rate given in \Cref{theorem: bundle complexity 1/epsilon log(epsilon)} and the analysis of Liang and Monteiro may be tighter than that of \Cref{{thm: 1/e convergence when g is smooth}}.

\Cref{{thm: 1/e convergence when g is smooth}} builds upon \cite{lacostejulien2015global}, which demonstrates a geometry-dependent rate of $\cO\left(D^2\gamma^{-2}(\log\epsilon^{-1})\right)$ for some FW variants on polyhedral sets. When the problem must be solved with high accuracy ($\epsilon\ll n^{-1}$), their result is more favorable than the previously known, dimension-independent, rate of $\cO({\epsilon}^{-1})$ for general closed convex sets. Likewise, our geometry-dependent rate is a substantial progress and reveals favorable convergence behavior of MPB-FA for seeking high-accuracy solutions for a broad class of convex composite nonsmooth problems.

The improved  {complexity} given in \Cref{thm: 1/e convergence when g is smooth} is possible at the cost  of including geometry-dependent constants that can be arbitrarily large and supposing that $f$ is piecewise linear and $g$ is smooth. From a dual perspective, this trade-off in dimension dependence for improved complexity in terms of $\epsilon$ is made possible by the structure of problem \eqref{problem: Dual}. It remains an open question whether a Frank-Wolfe-type algorithm with a similar trade-off can be designed when the feasible region is a general compact convex set.
More generally, problem \eqref{problem: Dual} is a strongly convex but nonsmooth problem for which none of the classical Frank-Wolfe variants  {are well suited}. Designing a Frank-Wolfe type method with an $\cO(\epsilon^{-1})$ convergence rate for problems like \eqref{problem: Dual} (nonsmooth, strongly-convex) was recently identified as an open problem in \cite{asgari2023projectionfree}.

\subsection*{Conclusions}
Expanding on the duality between Kelley's method and Fully-Corrective Frank-Wolfe to non-homogeneous convex piecewise linear functions, our work addresses a broader class of problems. Our interpretation of \textit{null steps} in a modified proximal bundle algorithm (MPB-FA) as dual to the Frank-Wolfe algorithm leads to a convergence rate of $\cO(\epsilon^{-1}(\log\epsilon^{-1}))$, surpassing the previous $\cO(\epsilon^{-2})$ rate. Further refinement achieves an iteration complexity of $\cO({\epsilon^{-8/9}}(\log\epsilon^{-1})^{2/9})$, a notable improvement. We also provide theoretical ground and preliminary numerical results to support the active-cut policy for bundle management. It is still an open question whether the $\cO({\epsilon^{-3}})$ rate of the PBM \cite{Kiwiel2000}
can be improved. We hope the investigation of modified bundle methods, such as MPB-FA, can ultimately lead to an improved complexity analysis of the classical PBM.

\appendix

\section{Technical proofs}

\subsection{Technical lemma}
\label{appendix: technical lemma}
The following lemma presents a result that, to the best of our knowledge, has not been explicitly stated in the existing literature. Specifically, we demonstrate  {a strong convexity type property for} the sum of a convex function and a strongly convex function defined on the Cartesian product of their domains.

\begin{lemma}
    \label{lm: Strongly convex in a variable}
    Let $p :{X} \subset \R^n\rightarrow \R$ be an $\alpha$-strongly convex differentiable function and $q: Y\subseteq \R^n\rightarrow \R$ be a convex differentiable function. If $(x^*, y^*)$ minimizes $p(x)+q(y)$ over $(x,y)\in Q\subseteq X\times Y$ a convex set, then
    \begin{align}
        \forall (x,y) \in Q, \quad  p(x) + q(y) \geq p(x^*) + q(y^*) + \frac{\alpha}{2}\|x-x^*\|^2.
    \end{align}
\end{lemma}
Remark: Note that the domains of $p$ and $q$ do not lie in the same space. So $p(x)+q(y)$ is not strongly convex jointly in $(x,y)$.
\begin{proof}
By the first order optimality condition, $\langle \nabla p(x^*), x-x^*\rangle + \langle \nabla q(y^*), y-y^*\rangle\geq 0$ for all $(x,y)\in Q$.
\begin{align*}
    p(x) + q(y) &= p(x^*) + \int_0^1{\nabla p(tx + (1-t)x^*)^\top (x-x^*)dt} + q(y)\\
    &\geq p(x^*) + \int_0^1{\langle \nabla p(tx + (1-t)x^*) - \nabla p(x^*), (x-x^*)\rangle dt} \\
    &\quad + q(y^*) + \nabla q(y^*)^\top  (y-y^*) + \nabla p(x^*)^\top (x-x^*)\\
    &\geq p(x^*) + q(y^*) + \int_0^1{\frac{1}{t}\langle \nabla p(tx + (1-t)x^*) - \nabla p(x^*), (tx-tx^*)\rangle dt}\\
    &\geq p(x^*) + q(y^*) + \int_0^1{\frac{1}{t}\alpha\|tx-tx^*\|^2 dt}= p(x^*) + q(y^*)  + \frac{\alpha}{2}\|x-x^*\|^2.
\end{align*}
The convexity of $q$ gives the first inequality, the second inequality follows from the optimality of $(x^*,y^*)$, and the third inequality uses the strong convexity of $p$.
\end{proof}
\subsection{Proof of \Cref{lemma: strong convexity constant of psi}}\label{appendix:lemma}
\begin{proof}[Proof of \Cref{lemma: strong convexity constant of psi}]
We will adapt the proof of Lemma 2.5 of \cite{beck2015linearly} to our setup. Let  $G = \max_{w \in E_1^{n}  {\conv(\cV)}}(\|\nabla g^*(-w)\|)$. 
 {Under \Cref{assum:g strong conv}, by \Cref{lemma:Primal-Dual variables correspondence},}
\begin{align}
    g^*(-w) &= -w^\top  x - g(x), \quad w = -\nabla g(x), \quad  {\text{and} \;\; x = \nabla g^*(-w)}. \label{eq:unique_w_x}
\end{align}
Using $\mu_g$-strong convexity of  {$g$}, 
\begin{align*}
    \|\nabla g(x) - \nabla g(x^*)\|  &\geq \mu_g \|x - x^*\| \\
    \frac{\|w^* - w\| }{\mu_g} &\geq \|\nabla g^*( {-w}) - x^*\| &&  {(\text{By \eqref{eq:unique_w_x}})}\\
    &\geq \|\nabla g^*( {-w})\| -\| x^*\|.
\end{align*}
After rearranging, we have $\|\nabla g^*(-w)\| \leq  {\frac{D_w }{\mu_g}} + \|x^*\|$,  {where $D_w$ is the diameter of \(\proj_w(\conv(\cV))\)}. As this inequality \textcolor{blue}{is} valid for any $w \in E_1^{n} \cV$, we have 
\begin{align*}
    G \leq  {\frac{D_w}{\mu_g}} + \|x^*\|. \stepcounter{equation}\tag{\theequation}\label{ineq: upper bound on G}
\end{align*}

 {By definition, $\psi(w,\beta) = g^*(-w)-\beta$, we have}
\begin{align}
    & \psi(w^*,\beta^*) - \psi(w,\beta)  - \langle \nabla \psi (w, \beta), ((w^*)^\top , \beta^*)^\top  - (w^\top , \beta)^\top  \rangle \notag \\
    = &  {\, g^*(-w^*)-g^*(-w)-\langle -\nabla g^*(-w), w^*-w\rangle} \notag  \\
    \ge & \frac{1}{2 L_g} \left\|w - w^*\right\|^2 \quad  {\text{(By $L_g^{-1}$-strong convexity of $g^*$)}}.
    \label{ineq : strong convexity in a variable for psi}
\end{align}
We bound the distance $|\beta^* - \beta|$ by the difference in function value $\psi(w,\beta) - \psi(w^*,\beta^*)$: 
\begin{align*}
    \beta - \beta^* &= \langle e_{n+1}, (w^\top , \beta)^\top  -  ((w^*)^\top , \beta^*)^\top \rangle\\
    &= \langle \nabla \psi(w^*, \beta^*), (w^\top , \beta)^\top  -  ((w^*)^\top , \beta^*)^\top \rangle - \langle \nabla g^*(-w^*), w - w^*\rangle \stepcounter{equation}\tag{\theequation}\label{eq: Beck upper bound on beta difference}\\
    &\leq \psi(w, \beta) - \psi(w^*, \beta^*) + \|\nabla g^*(-w^*)\| \|w - w^*\| \\
    &\leq  \psi(w, \beta) - \psi(w^*, \beta^*) + G \|w - w^*\|. \stepcounter{equation}\tag{\theequation}\label{eq: Beck upper bound on beta difference other side}
\end{align*}
Then, using \eqref{eq: Beck upper bound on beta difference} and Cauchy-Schwarz, we get $\beta^* - \beta \leq G \|w - w^*\|$, which combined with \eqref{eq: Beck upper bound on beta difference other side} gives
\begin{align}
   &\quad(\beta - \beta^*)^2\notag\\
   & {\le \max\{(G\|w-w^*\|)^2, (\psi(w, \beta) - \psi(w^*, \beta^*) + G \|w - w^*\|)^2\}}\nonumber\\
   &=  (\psi(w, \beta) - \psi(w^*, \beta^*) + G \|w - w^*\|)^2 \quad  {\text{(By \(\psi(w,\beta)-\psi(w^*,\beta^*)\ge 0\))}}\nonumber\\
   &=(\psi(w, \beta) - \psi(w^*, \beta^*))^2  {+} 2 G \|w {-} w^*\|(\psi(w, \beta) {-} \psi(w^*, \beta^*)) + G^2 \|w {-} w^*\|^2\nonumber\\
   & {\leq (\psi(w, \beta) - \psi(w^*, \beta^*)) (\psi(w, \beta) {-} \psi(w^*, \beta^*)+ 2 G D_w)} \label{ineq: geom strong conv 0}\\
   &\quad + G^2 2 L_g(\psi(w, \beta) - \psi(w^*, \beta^*))\nonumber\\
   & {=(\psi(w, \beta) - \psi(w^*, \beta^*)) (g^*(-w) - g^*(-w^*) + \beta^* - \beta + 2 G D_w)}\notag\\
   &\quad + G^2 2 L_g(\psi(w, \beta) - \psi(w^*, \beta^*))\nonumber\\
   &\leq (\psi(w, \beta) - \psi(w^*, \beta^*)) (\|\nabla g^*(-w) \|D_w + D_b+ 2 G D_w) \label{ineq: geom strong conv 1}\\
   &\quad + G^2 2 L_g(\psi(w, \beta) - \psi(w^*, \beta^*))\nonumber\\
   &\leq(\psi(w, \beta) - \psi(w^*, \beta^*)) (D_b + 3 G D_w  +2 L_g G^2 )\nonumber\\
   &\leq - \langle \nabla \psi (w, \beta), ((w^*)^\top , \beta^*)^\top  - (w^\top , \beta)^\top  \rangle (D_b + 3 G D_w  +2 L_g G^2).\label{ineq: geom strong conv 2}
\end{align}
Inequality \eqref{ineq: geom strong conv 0} uses \Cref{lm: Strongly convex in a variable}  {for \(\psi(w, \beta) - \psi(w^*, \beta^*)\ge (2L_g)^{-1}\|w-w^*\|^2\)}. For \eqref{ineq: geom strong conv 1}, we use the convexity of $g^*$. For 
\eqref{ineq: geom strong conv 2}, we use the convexity of $\psi$. Adding \eqref{ineq: geom strong conv 2} and \eqref{ineq : strong convexity in a variable for psi}, we get : 
\begin{align*}
   &\psi(w^*,\beta^*) -  \psi(w,\beta)  - 2 \langle \nabla \psi (w, \beta), ((w^*)^\top , \beta^*)^\top  - (w^\top , \beta)^\top  \rangle\\ &\quad\quad \geq (\max\{D_b + 3 G D_w  +2 L_g G^2, 2 L_g\})^{-1} (\|w - w^*\|^2 +(\beta - \beta^*)^2)\\
    &\quad\quad\geq (D_b + 3 G D_w  +2 L_g (G^2+1))^{-1}\left\|\binom{w^*}{ \beta^*} - \binom{w}{ \beta} \right\|^2. \stepcounter{equation}\tag{\theequation}\label{ineq: lower bound on geometric QG constant}
\end{align*}

This leads to the claimed inequality using \eqref{ineq: lower bound on geometric QG constant} and replacing $G$ by \eqref{ineq: upper bound on G} :
\begin{align*}
   &\psi(w^*,\beta^*) -  \psi(w,\beta)  - 2 \langle \nabla \psi (w, \beta), ((w^*)^\top , \beta^*)^\top  - (w^\top , \beta)^\top  \rangle\\ 
    &\geq \left( D_b + 3 \frac{4 M_f^2}{\mu_g}  +6M_f \|x^*\| + 2L_g \left[\left(\frac{2M_f}{\mu_g} 
    + \|x^*\|\right)^2 + 1\right]\right)^{-1}  \left\|\binom{w^*}{ \beta^*} - \binom{w}{ \beta} \right\|^2.
\end{align*}
\end{proof}

\bibliographystyle{siamplain}
\bibliography{references}
\end{document}

%% file: ex_shared.tex
\usepackage{lipsum}
\usepackage{amsfonts}
\usepackage{graphicx}
\usepackage{epstopdf}

\usepackage[utf8]{inputenc} 
\usepackage[T1]{fontenc}    
\usepackage{hyperref}       
\usepackage{url}            
\usepackage{booktabs}       
\usepackage{amsfonts}       
\usepackage{nicefrac}       
\usepackage{microtype}      
\usepackage{xcolor}         
\usepackage{apptools}
\usepackage{caption}
\usepackage{subcaption}
\usepackage{tikz, tikz-cd}
\usepackage{array, amsmath, mathtools, amssymb}
\usepackage{algorithm}
\usepackage{algpseudocode}
\usepackage{cite}
\usepackage{cancel}
\usepackage[normalem]{ulem}

\newcommand{\N}{\mathbb{N}}
\newcommand{\R}{\mathbb{R}}

\newcommand{\cV}{\mathcal{V}}
\newcommand{\cM}{\mathcal{M}}
\newcommand{\cO}{\mathcal{O}}

\newcommand{\conv}{\mathrm{conv}}

\newcommand{\proj}{\mathrm{Proj}}

\newcommand{\rank}{\mathrm{rank}}

\DeclareMathOperator*{\argmax}{arg\,max}
\DeclareMathOperator*{\argmin}{arg\,min}

\ifpdf
  \DeclareGraphicsExtensions{.eps,.pdf,.png,.jpg}
\else
  \DeclareGraphicsExtensions{.eps}
\fi

\newsiamremark{remark}{Remark}
\newsiamremark{hypothesis}{Hypothesis}
\crefname{hypothesis}{Hypothesis}{Hypotheses}
\newsiamthm{claim}{Claim}
\newsiamthm{assumption}{Assumption}

\title{An Absolute-Error Proximal Bundle Method through the Lens of Frank-Wolf}

\author{David Fersztand \thanks{Operations Research Center, MIT Cambridge, MA 
  }
\and Xu Andy Sun\thanks{Sloan School of Management and Operations Research Center, MIT, Cambridge, MA 
  (\email{sunx@mit.edu}).}
}


%% file: proof_final_improved_rate_v3.tex
\begin{proof}[Proof of \Cref{thm: 1/e convergence when g is smooth}]
The proof is made by case analysis.  {We first define the following constant $A$ by}
\begin{align}
    A^2 &:= \frac{\sqrt{2}\|x_0 - x^*\|\gamma}{6(1+5(\alpha\rho)^{-1})D\sqrt{\rho}} \epsilon^{-1/2}. \label{eq: optimal value of A}
\end{align}
 {It will become clear in Case 1 below why we define $A$ in this way. }

 {Now we analyze two cases: 1) $A\ge \sqrt{\alpha/2}$ and 2) $A<\sqrt{\alpha/2}$. In the first case, we can distinguish consecutive vs non-consecutive serious steps 
to carry out the blueprint \eqref{blueprint: improved_way} on Page \pageref{blueprint: improved_way}. In the second case, we do not distinguish the two types of serious steps and carry out the standard convergence analysis outlined in \eqref{blueprint: standard_way} on page \pageref{blueprint: standard_way}. Then, we combine both cases to reach the final improved convergence rate.}
\paragraph{\textbf{Case 1: $A\ge \sqrt{\alpha/2}$}}  {We draw a Venn diagram in \Cref{fig:venn}. The rectangle is the set of all serious steps. Circle $S_1$ contains all serious steps with ``large gradients'', i.e. \( \|y_k\| \geq \sqrt{\epsilon} M \), where $M:=\frac{6 A(1 + 5(\alpha \rho)^{-1}) D}{\gamma}$. Circle $S_2$ contains all serious steps with ``close centers'' \(\|x_{k-1} - x_{k}\| \le A \sqrt{\epsilon}\). Then, by \Cref{lm: next step is serious step}, every serious step in their intersection $S_3:=S_1\cap S_2$ is a consecutive serious step. 
}
\input{venn_case1}
 {Moreover, by \Cref{lm: number of serious steps with distant consecutive iterates}, the number of {serious steps} with ``distant centers''  \(\|x_{k-1} - x_{k}\| > A \sqrt{\epsilon}\) is at most \(\frac{2 \|x_0 - x^*\|^2}{\epsilon A^2} + \frac{2}{A^2\rho}\). By \Cref{lm: upper bound on the number of serious steps with small gradient}, the number of serious steps with ``small gradients'' is at most $\rho M^2$. If we use $|S^c|$ to denote the cardinality of the complement of a set $S$, then $|S_1^c|\le \rho M^2$, $|S_2^c|\le \frac{2 \|x_0 - x^*\|^2}{\epsilon A^2} + \frac{2}{A^2\rho}$. Note that the union $S_1^c \cup S_2^c$ contains all the \textit{non-consecutive} serious steps, which can be upper bounded by \(|S_1^c \cup S_2^c|\le \rho M^2 + \frac{2 \|x_0 - x^*\|^2}{\epsilon A^2} + \frac{2}{A^2\rho}\).}
We can further upper bound the number of \textit{consecutive} {serious steps} by the total number of {serious steps} $\frac{\rho \|x_0 - x^*\|}{\epsilon} + 1$ given in \Cref{lm: Number of serious steps}. We use the upper bound on the number of \textit{null steps} between two \textit{serious steps} given by \Cref{lemma: number of null steps}.  {Then, summarize all the above counting, we can implement the blueprint \eqref{blueprint: improved_way} and get the following upper bound on the total number of steps of \Cref{alg: Proximal Bundle Method}:}
\begin{align}
\underbrace{1+\frac{\rho \|x_0 - x^*\|^2}{\epsilon}}_{\text{consecutive serious steps}} +& \left(\hspace{-0.1cm}1{+}\underbrace{\rho \hspace{-1.6pt}\left(\frac{6(1{+}5(\alpha\rho)^{-1}) AD}{\gamma}\right)^2 \hspace{-0.25cm}}_{\text{after small gradient occurs}}{+}\underbrace{\frac{2 \|x_0 {-} x^*\|^2}{\epsilon A^2} {+} \frac{2}{A^2\rho}}_{\text{distant center steps}}\hspace{-0.1cm}\right)\times\label{eq:total_iter_case1}
    \\
& \hspace{5mm}\underbrace{\left[1{+}\max\left\{2,\frac{D^2}{\bar{\mu}_{\psi, \rho} \rho \gamma^2}\right\} \log \left(\frac{4D^4}{\epsilon^2\rho^2} \right)\right]}_{\text{number of consecutive null steps}}.\notag
\end{align}
 {A simple calculation verifies the value of $A$  in \eqref{eq: optimal value of A} minimizes} $\rho \left(\frac{6(1+5(\alpha\rho)^{-1}) AD}{\gamma}\right)^2 + \frac{2\|x_0 - x^*\|^2}{\epsilon A^2}$.  {This is why we choose $A$ according to \eqref{eq: optimal value of A}.}

\paragraph{\textbf{Case 2: $A < \sqrt{\alpha/2}$}} 
Following \eqref{eq: optimal value of A}, 
\begin{align*}
    \frac{\sqrt{2}\|x_0 - x^*\|\gamma}{6(1+5(\alpha\rho)^{-1})D\sqrt{\rho}} \epsilon^{-1/2}&< \frac{ \alpha}{2}.
\end{align*}
This gives an inequality involving the upper bound of the number of \textit{serious steps}. 
\begin{align*}
    \frac{\rho \|x_0 - x^*\|^2}{\epsilon}&<\frac{6^2(1+5(\alpha\rho)^{-1})^2 \alpha^2 D^2 \rho^2}{2\cdot 4 \cdot  \gamma^2}
    =\frac{9(\alpha\rho+5)^2D^2 }{2\gamma^2}\\
    &\leq \frac{9(1 + \frac{\rho}{2L_g} + 5)^2D^2 }{2\gamma^2} = \frac{D^2}{\gamma^2}  \left(162 +  {\frac{27\rho}{L_g}}+\frac{9\rho^2}{ {8}L_g^2}\right).
\end{align*}
In this case,  {$S_3$ may contain both consecutive and non-consecutive serious steps. We need to use the standard way \eqref{blueprint: standard_way}} to upper bound the total number of iterations by
\begin{align}
    &\left(\frac{\rho \|x_0 - x^*\|^2}{\epsilon} + 1\right)\left[1+\max\left\{2,\frac{D^2}{\bar{\mu}_{\psi, \rho} \rho \gamma^2}\right\} \log \left(\frac{4D^4}{\epsilon^2\rho^2} \right)\right]\notag\\
    &\quad\quad\quad\leq \frac{D^2}{\gamma^2} \left(162 +  {\frac{27\rho}{L_g}}+\frac{9\rho^2}{ {8}L_g^2}\right)\left[1+\max\left\{2,\frac{D^2}{\bar{\mu}_{\psi, \rho} \rho \gamma^2}\right\} \log \left(\frac{4D^4}{\epsilon^2\rho^2} \right)\right].\label{eq:total_iter_case2}
\end{align}
Putting these two bounds \eqref{eq:total_iter_case1} and \eqref{eq:total_iter_case2} together, we obtain an upper bound on the total number of iterations that is independent of $A$,
\begin{align}
    &\underbrace{1 +\frac{\rho \|x_0 - x^*\|^2}{\epsilon}}_{\text{consecutive serious steps}}  + \underbrace{\left[1+\frac{K\sqrt{\epsilon}}{2\|x_0 - x^*\|^2\sqrt{\rho}}+\frac{K\sqrt{\rho}}{\sqrt{\epsilon}} {+} \frac{D^2}{\gamma^2}\left(162 +  {\frac{27\rho}{L_g}}+\frac{9\rho^2}{ {8}L_g^2}\right)\right]}_{\text{serious steps with {null steps} inbetween}}\times\label{final convergence rate last thm}\\
    &\hspace{3.5cm}  \underbrace{\left[1{+}\max\left\{2,\frac{D^2}{\bar{\mu}_{\psi, \rho} \rho \gamma^2}\right\} \log \left(\frac{4D^4}{\epsilon^2\rho^2} \right)\right]}_{\text{number of consecutive null steps}}, \nonumber
\end{align}
where $K:= 2\cdot\frac{6\sqrt{2}(1+5(\alpha\rho)^{-1})D\|x_0 - x^*\|}{\gamma }$.  {Since 
$K$ contains $(\alpha\rho)^{-1}$ and $\alpha$ is a function of $\rho$ by \Cref{lm: strong convexity of phi plus proximal}, let us separate $(\alpha\rho)^{-1}$ out and write $K$ as
\begin{align}
    K = K_1 + K_2 \frac{1}{\alpha\rho},
\end{align}
where $K_1, K_2$ contains all the parameters other than $\epsilon, \rho, \alpha$. By Lemma 3.16, we have
\begin{align}
2{\alpha\rho} = (\rho/L_g+2)-\sqrt{\rho^2/L_g^2+4} = \frac{4\rho/L_g}{\rho/L_g + 2 +\sqrt{\rho^2/L_g^2 + 4}}, 
\end{align}
Therefore, $\frac{1}{\alpha\rho} = \frac{\rho/L_g + 2 +\sqrt{\rho^2/L_g^2 + 4}}{2\rho/L_g} = \mathcal{O}(\frac{1}{\rho})$, as $\rho\to 0$.
Then, $K$ can be treated as $K = K_1 + \tilde{K}_2 \frac{1}{\rho}$,
where $\tilde{K}_2$ is independent of $\alpha$. 
Also from (3.19), we know $\bar{\mu}_{\psi,\rho}=\mathcal{O}(\frac{1}{\rho^2})$ as $\rho\to 0$. 
Substituting $K$ and $\bar{\mu}_{\psi,\rho}$ into 
\eqref{final convergence rate last thm}, we get
\begin{align}
    &C_0 \frac{\rho}{\epsilon} + \biggl(C_1 \frac{\sqrt{\epsilon}}{\sqrt{\rho}} + C_2\frac{\sqrt{\epsilon}}{\rho^{3/2}} + C_1 \frac{\sqrt{\rho}}{\sqrt{\epsilon}} + C_2\frac{1}{\sqrt{\epsilon\rho}} + C_3(\rho) \biggr)\cdot\biggl(\frac{1}{\rho^3} \log(\frac{1}{\epsilon\rho})\biggr),\label{eq:upperbound_order}
\end{align}
where $C_i$ for $i=0,1,2$ are constants independent of $\epsilon$ and $\rho$, and $C_3(\rho)$ is the quadratic function of $\rho$ for the term $\frac{D^2}{\gamma^2}\left(162 +  {\frac{27\rho}{L_g}}+\frac{9\rho^2}{ {8}L_g^2}\right)$.}
 {
Since $\log(1/\epsilon)$ and $\log(1/\rho)$ are on the same order, we can replace $\log(1/\epsilon\rho)$ by $\log(1/\epsilon)$. Taking only the highest order term $\frac{1}{\sqrt{\epsilon\rho}}$ in the regime of $\rho\to 0$ and $\epsilon\to 0$, we can write \eqref{eq:upperbound_order} as
\begin{align}
    C_0 \frac{\rho}{\epsilon} + C_1 \frac{\rho^{-7/2}}{\sqrt{\epsilon}}\log(\frac{1}{\epsilon}),\label{eq:asympt1}
\end{align}
which is a convex function in $\rho$. Taking the minimum by setting the derivative to zero, we have $C_0 \frac{1}{\epsilon} = \frac{7C_2}{2} \frac{\rho^{-9/2}}{\sqrt{\epsilon}}\log(\frac{1}{\epsilon})$, which implies $\rho = \mathcal{O}\biggl(\epsilon^{1/9}\log({\epsilon}^{-1})^{2/9}\biggr)$.
Substituting back to \eqref{eq:asympt1}, the upper bound for the total number of iterations is on the order of $\mathcal{O}(\epsilon^{-8/9}\log({\epsilon}^{-1})^{2/9})$.}
\end{proof}

%% file: venn_case1.tex
\def\firstcircle{(0,-0.5cm) circle (2.5cm)}
\def\thirdcircle{(3.5cm,-0.5cm) circle (2.5cm)}

\begin{figure}[htb]
    \centering
\begin{tikzpicture}[scale=0.6, every node/.style={scale=0.6}]
    \draw (-3cm, -3.5cm) rectangle (6.5cm, 3cm);
    \draw (1.5cm, 3cm) node [above] {all serious steps };
    \draw \firstcircle node [below] {$S_1$};
    \draw (-0.5cm,0cm) node {(large gradient)};
    \draw (-0.5cm, 0.5cm) node {$y_k\ge\sqrt{\epsilon}M$};
    \draw (4cm, 0.5cm) node {$\|x_{k-1}-x_k\|\le A\sqrt{\epsilon}$};
    \draw (4cm, 0cm) node {(close centers)};
    \draw \thirdcircle node [below] {$S_2$};


    \begin{scope}
      \clip \firstcircle;
      \fill[blue!10] \thirdcircle;
      \draw (1.8cm, -0.5cm) node [below] {$S_3$};
    \end{scope}
\end{tikzpicture}
\caption{{The Venn diagram for the set of serious steps.} }
    \label{fig:venn}
\end{figure}